\documentclass[12pt,a4paper]{article}
\usepackage[utf8]{inputenc}
\usepackage{xifthen}
\usepackage{csquotes}
\usepackage{amssymb}
\usepackage{amsthm}
\usepackage{mathtools}
\usepackage{physics} 
\usepackage{enumerate}
\usepackage{hyperref}
\usepackage{cleveref}
\usepackage{nameref}
\usepackage{comment}
\usepackage{times}
\newcommand{\elminus}[1]{\setminus\set*{#1}}

\DeclareMathAlphabet\mathbfcal{OMS}{cmsy}{b}{n}

\DeclarePairedDelimiterX\set[1]\lbrace\rbrace{\let\given\undefined\newcommand\given{\;\delimsize\vert\;}#1}
\DeclarePairedDelimiter{\prn}{\lparen}{\rparen}

\newcommand{\tidx}[1]{{(#1)}}

\newcommand{\N}{\mathbb{N}}
\newcommand{\R}{\mathbb{R}}
\newcommand{\collc}{\mathbfcal{C}}

\newcommand{\collq}{\mathbfcal{Q}}

\newcommand{\diff}{\mathop{}\!\mathrm{d}}

\newcommand{\ball}[1]{{B^{#1}}}
\newcommand{\stdb}{\ball n}

\newcommand{\sphere}[1]{\mathbb{S}^{#1}}
\newcommand{\stdsph}{\sphere{n-1}}

\DeclareMathOperator{\im}{im}

\DeclareMathOperator{\cl}{cl}
\DeclareMathOperator{\V}{V}
\DeclareMathOperator{\Su}{S}

\newcommand{\SufC}{\Su_{f, \collc}}
\DeclareMathOperator{\conv}{conv}
\DeclareMathOperator{\mPos}{pos}
\DeclareMathOperator{\diam}{diam}
\DeclareMathOperator{\vspan}{span}
\DeclareMathOperator{\pspan}{\overline{span}}
\DeclareMathOperator{\supp}{supp}

\newcommand{\K}{\mathcal{K}}
\newcommand{\Pt}{\mathcal{P}}
\newcommand{\mbS}{\mathbb{S}}

\newtheorem{baseenv}{Abstract Env}[section]
\newtheorem{theorem}[baseenv]{Theorem}
\newtheorem{definition}[baseenv]{Definition}
\newtheorem{remark}[baseenv]{Remark}
\newtheorem{lemma}[baseenv]{Lemma}
\newtheorem{corollary}[baseenv]{Corollary}
\newtheorem{example}[baseenv]{Example}
\newtheorem{proposition}[baseenv]{Proposition}
\newtheorem*{theorem*}{Theorem}

\title{Extremizers of the Alexandrov--Fenchel inequality within a new class of convex bodies}
\author{Daniel Hug and Paul A. Reichert}
\date{\today}

\begin{document}

\maketitle

\begin{abstract}
    Mixed volumes in $n$-dimensional Euclidean space are functionals of $n$-tuples consisting of convex bodies $K,L,C_1,\ldots,C_{n-2}$. The Alexandrov--Fenchel inequalities are fundamental inequalities between mixed volumes of convex bodies, which cover as very special cases many important inequalities between basic geometric functionals. The problem of characterizing completely the equality cases in the Alexandrov--Fenchel inequality is wide open. Major recent progress was made by Yair Shenfeld and Ramon van Handel   \cite{SvH22,SvH23+}, in particular they resolved the problem in the cases where $K,L$ are general convex bodies and $C_1,\ldots,C_{n-2}$ are polytopes, zonoids or smooth bodies (under some dimensional restriction).
    We introduce the class of polyoids, which includes polytopes, zonoids and triangle bodies, and characterize polyoids by  using generating measures.
    Based on this characterization and Shenfeld and van Handel's contribution, we extend their result to a class of convex bodies containing all polyoids and smooth bodies. Our result is stated in terms of the support of the mixed area measure of the unit ball $B^n$ and $C_1,\ldots,C_{n-2}$. A geometric description of this support is provided in the accompanying work \cite{HugReichert23+}.
\end{abstract}

\paragraph{\small MSC-classes 2020.}{\small 52A39, 52A20, 52A21 52A40}

\section{Introduction}
Mixed volumes of convex bodies (nonempty compact convex sets) in Euclidean space $\R^n$, $n\ge 2$, are symmetric functionals of $n$-tuples of convex bodies. They naturally arise
as coefficients of polynomial expansions of nonnegative Minkowski combinations of convex bodies. Writing $\V$ for the volume functional (Lebesgue measure) and $\alpha_1K_1+\cdots+\alpha_mK_m$ for the Minkowski combination of the convex bodies $K_1,\ldots,K_m\subset\R^n$ with nonnegative coefficients $\alpha_1,\ldots,\alpha_m\in\R$, we have
\begin{equation}\label{eq:1.1}
\V(\alpha_1K_1+\cdots+\alpha_mK_m)=\sum_{i_1,\ldots,i_n=1}^m \V(K_{i_1},\ldots,K_{i_n})\alpha_{i_1}\cdots\alpha_{i_n},
\end{equation}
where $\V(K_{i_1},\ldots,K_{i_n})$ is called the mixed volume of $K_{i_1},\ldots,K_{i_n}$. As symmetric functions of their $n$ arguments, mixed volumes are uniquely determined by this expansion. We refer to \cite[Chap.~5.1]{Schneider} or
\cite[Chap.~3.3]{Hug} for an introduction to mixed volumes. Conversely, the mixed volume $\V(K_1,\ldots,K_n)$ of a given $n$-tuple of convex bodies $K_1,\ldots,K_n$ can be obtained as an alternating sum of  volumes of Minkowski sums, that is,
\begin{equation}\label{eq:1.2}
\V(K_1,\ldots,K_n)=\frac{1}{n!}\sum_{k=1}^n(-1)^{n+k}\sum_{1\le i_1<\cdots<i_k\le n}\V(K_{i_1}+\cdots +K_{i_k}).
\end{equation}
While relations \eqref{eq:1.1} and \eqref{eq:1.2} can be efficiently employed for introducing mixed volumes and understanding some of their basic properties, their usefulness in deriving inequalities for mixed volumes seems to be limited. We refer to Schneider \cite[Notes for
Sect.~5.1]{Schneider} for further background information. A deep inequality for mixed volumes of convex bodies, with many consequences and applications to diverse fields, has been found and established by Alexandrov \cite{AF1937} (see Schneider \cite[Notes for Sect. 7.3]{Schneider}, also for some historical comments).

\begin{theorem*}[Alexandrov--Fenchel Inequality]
    \makeatletter\def\currentlabelname{Alexandrov--Fenchel Inequality}\makeatother\label{thm:af}
    Let  $K, L\subset\R^n$ be convex bodies, and let $\collc = (C_1, \ldots, C_{n-2})$ be an $(n-2)$-tuple of convex bodies in $\R^n$. Then
    \begin{align}
      \V(K, L, \collc)^2 \ge \V(K, K, \collc) \V(L, L, \collc),\tag{AFI}
    \end{align}
where $\V(K, L, \collc):=\V(K,L,C_1, \ldots, C_{n-2})$.
\end{theorem*}

We state the Alexandrov--Fenchel inequality in a second version. It makes use of a linear extension of mixed volumes,  known already to Alexandrov, to differences of support functions of convex bodies (see \cite[Sect.~5.2]{Schneider}). Such extensions turned out to be useful in  proofs of the inequality. For a convex body $K\subset\R^n$, the support function $h_K:\R^n\to\R$ is defined by $h_K(u):=\max\{\langle x,u\rangle:x\in K\}$, where $\langle \cdot\,,\cdot\rangle$ denotes the Euclidean scalar product. The support function is positively homogeneous of degree one, hence it is often sufficient to consider its restriction to the Euclidean unit sphere $\mathbb{S}^{n-1}$, and it uniquely determines the underlying convex body $K$.

\begin{theorem*}[General Alexandrov--Fenchel Inequality]
    \makeatletter\def\currentlabelname{General Alexandrov--Fenchel Inequality}\makeatother\label{thm:af2}
    Let $K_1, K_2, L\subset\R^n$ be convex bodies, and let $\collc = (C_1, \ldots, C_{n-2})$ be an $(n-2)$-tuple of convex bodies in $\R^n$. Then
    \begin{align}
      \V(h_{K_1} - h_{K_2}, L, \collc)^2 \ge \V(h_{K_1} - h_{K_2}, h_{K_1} - h_{K_2}, \collc) \V(L, L, \collc).\tag{GAFI}
    \end{align}
\end{theorem*}

For a proof of the equivalence of the two versions, we refer to Shenfeld and van Handel \cite[Lem.~3.11]{SvH23+}.

\medskip

Despite considerable effort, to date  it is unknown when exactly equality holds in (AFI) for a general $(n-2)$-tuple of convex bodies $\collc = (C_1, \ldots, C_{n-2})$ in $\R^n$. While the equality cases can be described purely by dimensionality considerations when at least one of the mixed volumes vanishes, the situation turns out to be much more subtle when $\V(K, K, \collc)$ and $\V(L, L, \collc)$ are both positive.

Shenfeld and van Handel \cite{SvH23+} recently fully characterized the equality cases in (AFI) when $C_1, \ldots, C_{n-2}$ are polytopes. Then they used this result to achieve a characterization when  $C_1, \ldots, C_{n-2}$ are polytopes, zonoids or smooth convex bodies, under a mild dimensionality assumption, which they called supercriticality (see \cite[Thm.~14.9]{SvH23+}). Supercriticality  is a natural condition that provides some dimensional restriction on a given sequence of nonempty sets, which is satisfied e.g. for any sequence $\collc=(C_1,\ldots,C_{n-2})$ of  full-dimensional convex bodies in $\R^n$. It is related to a well-known  condition that ensures that the mixed volume of a given tuple of convex bodies is positive. We refer to Section \ref{sec:3} for a precise definition and basic properties that are related to this concept. Recall that a zonoid is a limit (with respect to the Hausdorff metric) of a sequence of finite Minkowski sums of segments. A convex body $K$ is said to be smooth if each boundary point of $K$ is contained in a unique supporting hyperplane of $K$.

In this article, we study a class of convex bodies, which we call
\emph{polyoids}, that encompasses  polytopes, zonoids and triangle bodies.  Polyoids are obtained as limits of sequences of polytopes  that are finite Minkowski sums of polytopes  having at most a fixed number $k$ of vertices (for some $k\in\N$). If $k=2$, we are back in the zonoid case. For $k=3$ we cover the class of triangle bodies (cf.~\cite[Sect. 3]{Schneider1996}, \cite[p.~201]{Schneider}). If $\Pt^n_k$ denotes the set of $k$-topes in $\R^n$ (polytopes having at most $k$ vertices), then the class of polyoids in $\R^n$ is the union of the Minkowski classes $\mathfrak{M}(\Pt^n_k)$, $k\in\N$,  generated by  $\Pt^n_k$ (see \cite[Sect.~3.5]{Schneider} for  information on Minkowski classes and additive generation). Our treatment will be limited to supercritical $(n-2)$-tuples of convex bodies in $\R^n$. We refer to Section \ref{sec:2} for  precise definitions and further discussion for these notions.  The main aim of this work is to extend the characterization of the equality cases for (AFI), obtained by Shenfeld and van Handel \cite{SvH23+}, to all convex bodies $K,L$ and all supercritical tuples $\collc = (C_1, \ldots, C_{n-2})$ of polyoids and smooth convex bodies.

We begin by stating Shenfeld and van Handel's result for supercritical tuples of polytopes, zonoids and smooth convex bodies. For this purpose, we need the mixed area measure $\Su(K_1,\ldots,K_{n-1},\cdot)$ of an $(n-1)$-tuple of convex bodies $K_1,\ldots,K_{n-1}\subset\R^n$. Recall from \cite[Sect.~5.1]{Schneider} (or \cite[Thm.~4.1]{Hug}) how these finite Borel measures on the Euclidean unit sphere $\mbS^{n-1}$ are defined. They are related to mixed volumes and support functions via the relation
\begin{equation}\label{eqex}
\V(K_1,\ldots,K_{n-1},K_n)=\frac{1}{n}\int_{\mathbb{S}^{n-1}} h_{K_n}(u)\, \Su(K_1,\ldots,K_{n-1},\diff u),
\end{equation}
which holds for all convex bodies $K_1,\ldots,K_n\subset \R^n$. For given $K_1,\ldots,K_{n-1}$, the mixed area measure $\Su(K_1,\ldots,K_{n-1},\cdot)$ is the unique Borel measure on $\mathbb{S}^{n-1}$ such that \eqref{eqex} holds for all convex bodies $K_n$.  As in the case of the mixed volume, also the mixed area measure can be extended as an $(n-1)$-linear map to differences of support functions (see again \cite[Sect.~5.2]{Schneider}). Then relation \eqref{eqex}
remains true with convex bodies   replaced by differences of support functions.
 If $B^n$ is the Euclidean unit ball, then $\supp\Su(\stdb, \collc,\cdot)$ denotes the support of the mixed area measure of $B^n$ and $\collc = (C_1,\ldots, C_{n-2})$, which is the complement of the largest open subset of $\mathbb{S}^{n-1}$ on which $\Su(\stdb, \collc,\cdot)$ vanishes.

\begin{theorem*}[\mbox{\cite[Thm. 14.9]{SvH23+}}]\label{thm:afcharPoly}
    Let $K, L \subset\R^n$ be convex bodies, and let $\collc = (C_1, \ldots, C_{n-2})$ be a supercritical $(n-2)$-tuple of polytopes, zonoids or  smooth convex bodies in $\R^n$ such that
    $\V(K, K, \collc), \V(L, L, \collc) > 0$.
    Then  {\rm (AFI)} holds with equality if and only if there are $a > 0$ and $x \in \R^n$ such that
    $h_K = h_{aL + x}$  on $\supp\Su(\stdb, \collc,\cdot)$.
\end{theorem*}

In the case where $C_1, \ldots, C_{n-2}$ are all smooth, the result was already known (see \cite[Thm.~7.6.8]{Schneider} and the comment after \cite[Thm.~1.2]{HugReichert23+}). The main point of the preceding theorem is that a mixture of smooth and non-smooth bodies (which then are polytopes or zonoids) is admitted.
 Shenfeld and van Handel also characterized the much more involved  equality cases for arbitrary tuples of polytopes $\collc$. For their treatment of the zonoid case, the characterization theorem for polytopes is used as a crucial ingredient.

 Our main result is an extension of the preceding theorem where zonoids and polytopes are included in the larger class of
 polyoids.

\begin{theorem}
\label{corfin}
Let $K, L \subset\R^n$ be convex bodies,  and let $\collc = (C_1, \ldots, C_{n-2})$ be a supercritical $(n-2)$-tuple of polyoids or  smooth convex bodies in $\R^n$.
\begin{enumerate}[{\rm (a)}]
\item If \ $\V(K,L,\collc)=0$, then {\rm (AFI)} holds with equality and $K,L$ are homothetic.
\item Let $\V(K,L,\collc)>0$. Then {\rm (AFI)} holds with equality if and only if there are $a>0$ and $x\in\R^n$ such that $h_K=h_{aL+x}$ on $
\supp \Su(\stdb,\collc,\cdot)$.
\end{enumerate}
\end{theorem}

At the end of  Section \ref{sec:2}, we introduce a formally larger class of convex bodies (which we call macroids) for which the statement of Theorem \ref{corfin} remains true.

In \cite[Sect.~4]{Schneider1988}, Schneider established a characterization of the equality cases in the Alexandrov--Fenchel inequality for convex bodies $K,L$ and zonoids $C_1,\ldots,C_{n-2}$, under the additional assumption that $K,L$ are centrally symmetric and all bodies are full-dimensional.
Schneider's characterization involves specific geometric information about $\collc$, namely the closure of the set of all extremal normal vectors of  the $(n-1)$-tuple $(\stdb, \collc)$. In contrast,  Shenfeld and van Handel first characterize the equality cases in terms
 of the support of the mixed area measure $\supp\Su(\stdb, \collc,\cdot)$ (without the assumption of central symmetry of $K,L$ and with a relaxed  dimensionality assumption). Finally, they show  that $\supp\Su(\stdb, \collc,\cdot)$ equals the closure of the set of all extremal normal vectors of the $(n-1)$-tuple  $(\stdb, \collc)$ (see \cite[Prop.~14.13]{SvH23+}).

 According to a general conjecture due to Schneider \cite[Conjecture 7.6.14]{Schneider},  for an arbitrary $(n-1)$-tuple of convex bodies $(C,\collc)$ the support of the mixed area measure $\supp\Su(C, \collc,\cdot)$ is precisely the closure of the set of all extremal normal vectors of $(C, \collc)$. This conjecture is open even in the case where all bodies are zonoids (see \cite[Sect.~4]{Schneider1988} for some discussion). For the application to the equality cases of the Alexandrov--Fenchel inequality, only the special case where $C=\stdb$ is required. In \cite{HugReichert23+}, Schneider's conjecture concerning the support of mixed area measures is settled for the class of polyoids, which in particular covers the case where all $n-1$ bodies are general zonoids. In combination with the results of the present work, we thus obtain a geometric characterization of the equality cases in (AFI) not only for  general convex bodies $K,L$ and zonoids, but for general convex bodies $K,L$ and the larger class of polyoids (and smooth bodies).

\medskip

The paper is structured as follows. In Section \ref{sec:2} we deduce a representation result for the support functions of polyoids, stated as Corollary  \ref{thm:ptb-char}, from a more general result concerning Minkowski classes generated by homothety invariant closed families of convex bodies.  A related representation theorem for support functions of zonoids in terms of their generating measures  is a well-known and versatile tool in the study of zonoids. We also introduce the larger class of macroids whose definition is motivated by Corollary \ref{thm:ptb-char}. In Section \ref{sec:3} we start with a brief discussion of supercritical tuples of sets. Then we prepare the proof of Theorem \ref{thm:supercritical}, which is an equivalent version of Theorem \ref{corfin}, that involves the mixed area measure of a difference of support functions. Our arguments are inspired by and  partly based on the results by Shenfeld and van Handel \cite{SvH23+}. Theorems \ref{thm:supercritical} and \ref{corfin} both hold within the formally  larger class of macroids. In Appendix \ref{app:macroid-not-polyoid} we construct a macroid that is not a polyoid.

\section{Polyoids and beyond}\label{sec:2}

In this section, we introduce the class of polyoids and establish a characterization theorem. Our definition is guided by the geometric definition of a zonoid as a limit of a sequence of zonotopes, where a zonotope is a finite Minkowski sum of segments.

In the following, we work in Euclidean space $\R^n$ with scalar product $\langle\cdot\,,\cdot \rangle$ and norm $\|\cdot\|$. For a set $A\subseteq\R^n$, we set $A^\perp:=\{x\in\R^n\colon \langle x,a\rangle=0 \text{ for }a\in A\}$, the linear subspace orthogonal to the linear span  of $A$, and $u^\perp:=\{u\}^\perp$ for $u\in\R^n$. The volume of the Euclidean unit ball $B^n$ is denoted by $\kappa_n$, its surface area is $\omega_n:=n\kappa_n$. If we write $\mathcal{H}^{n-1}$ for the $(n-1)$-dimensional Hausdorff measure in $\R^n$ and $\mathbb{S}^{n-1}$ for the unit sphere, then $\mathcal{H}^{n-1}(\mathbb{S}^{n-1})=\omega_n$ for $n\ge 1$.  Most of the time we focus on $n\ge 2$, but almost all  statements and definitions hold for $n\in\N_0$ (if properly interpreted).  We write $\K^n$ for the set of nonempty compact convex sets in $\R^n$ and endow $\K^n$  with the Hausdorff metric. Elements of $\K^n$ are called convex bodies.

A map $\varphi:\R^n\to\R^n$ is a dilatation if there is some $\lambda> 0$ such that $\varphi(x)=\lambda x$ for $x\in\R^n$. A homothety is a dilatation followed by a translation.

For $k\in\N$ we set $[k]:=\{1,\ldots,k\}$.
If $(E,\rho)$ is a metric space and $A\subseteq E$ is nonempty, then $\diam A:=\sup\{\rho(x,y) \;\colon x,y\in A\}\in [0,\infty]$ denotes the diameter of $A$.

\begin{definition}
 {\rm For each $k\in\N$, let $\Pt^n_k\subset\K^n$ be the set of polytopes in $\R^n$ with at most $k$ vertices. Elements of $\Pt^n_k$ are called \emph{$k$-topes}. A finite Minkowski sum of $k$-topes is called a \emph{$k$-polyotope}.}
\end{definition}

\begin{remark}{\rm \label{thm:ktopeClosed}
 For any compact set $A\subset\R^n$,
  the set $\set*{ P \in \Pt^n_k \;\colon P \subseteq A }\subset\K^n$
  is  compact. Hence $\Pt^n_k$ is a countable union of compact subsets of $\K^n$ and thus a  measurable subset of $\K^n$. It is convenient to consider the subspace $\sigma$-algebra on $\Pt^n_k$ which is induced by the Borel $\sigma$-algebra of $\K^n$.}
\end{remark}

Next we define a class of convex bodies which generalizes the class of zonoids and contains arbitrary polytopes.

\begin{definition}
\label{def:polyoid}
  {\rm Let $k \in \N$ and $K \in \K^n$.
  If $K\in\K^n$ is the limit of a sequence of $k$-polyotopes, then $K$ is called a \emph{$k$-polyoid}. A convex body  $K$ is called a \emph{polyoid} (a \emph{polyotope}) if it is a $k$-polyoid (a $k$-polyotope) for some $k\in\N$.}
\end{definition}

\begin{remark} {\rm \begin{enumerate}[{\rm (a)}]
\item For a given $k\in\N$, the class of $k$-polyoids in $\R^n$ is a closed subset of $\K^n$. In the terminology of \cite[Sect.~3.5]{Schneider}, the class of $k$-polyoids is the Minkowski class $\mathfrak{M}(\Pt^n_k)$ generated by $\Pt^n_k$.
\item  A $1$-polyoid is just a singleton, a $2$-polyoid is a zonoid and a $3$-polyoid is a \emph{triangle body}, as defined in \cite[p.~201]{Schneider} (or \cite[Sect.~3]{Schneider1996}). Moreover, for a given polytope $P$ there is some integer $k\in\N$ (depending on $P$) such that $P$ is a $k$-polyotope and hence a $k$-polyoid.
\item Clearly, $\Pt^n_k\subseteq \Pt^n_{\ell}$ for $k\le\ell$. Hence any $k$-polyoid is  an $\ell$-polyoid for $k\le \ell$. In particular, if $C_1,\ldots,C_{r}$ are polyoids in $\R^n$, for a fixed $r\in\N$, then there is some $k\in\N$ such that $C_1,\ldots,C_{r}$ are $k$-polyoids. Similar statements hold for polyotopes.
\item In $\R^2$ every centrally symmetric convex body is a 2-polyoid (a zonoid), and every convex body in $\R^2$ is a $3$-polyoid (a triangle body). The first fact is well-known (cf.~\cite[Cor. 3.5.7]{Schneider}), the second follows from \cite[Thm.~3.2.14]{Schneider}.
\item Let $n\ge 3$. If $k\in\N$ is fixed and $P_k^*$ is an indecomposable polytope in $\R^n$ with more than $k$ vertices, then it follows from \cite[Thm.~3.4.2]{Schneider} (see also \cite[Thm. 4]{Berg69}) that $P_k^*$ is not approximable by the class $\Pt^n_k$. Hence $P_k^*$ is not a $k$-polyoid, but certainly $P_k^*$ is an $\ell$-polyoid, for some $\ell>k$. For instance, for each $k\ge 2$, there is some indecomposable $(k+1)$-tope (with triangular $2$-faces) which is not a $k$-polyoid.
\item The Minkowski sum of a triangle in $\R^2\times\{0\}$ and a $2$-dimensional ball  in $\{0\}\times \R^2$ yields an example of a $3$-polyoid  which is not a zonoid, not a polytope, and  neither smooth nor strictly convex. It is clear from \cite[Cor.~3.5.12]{Schneider} that the class of $3$-polyoids is much larger than the class of zonoids.
\item For a given $k\in [n]$, Ricker \cite{Ricker81} calls a finite Minkowski sum of $r$-dimensional simplices with $r\in\{0,\ldots,k\}$ a $k$-zonotope. Each such $k$-zonotope is a particular $(k+1)$-polyotope, for $k\in[n]$. Ricker then defines a $k$-zonoid (for $k\in[n]$) as a limit of $k$-zonotopes and characterizes $k$-zonoids in terms of the ranges of $k$ vector measures, thus extending a known result for 1-zonoids (i.e., zonoids).  A $3$-dimensional double pyramid over a triangle base  is not a $k$-zonoid (as follows from \cite[Thm.~3.4.2]{Schneider}), for any $k\in \N$, but it is  a $5$-polyotope.
\item Let $K$ be an $n$-dimensional convex cone which is not a polytope.  Then $K$ is indecomposable by \cite[Thm.~2]{Sallee74}, hence  \cite[Thm.~3.4.2]{Schneider}  implies that $K$ is not a polyoid.
    \end{enumerate}
}
\end{remark}

We now prepare the proof of Corollary \ref{thm:ptb-char}, which is an analogue for polyoids of a well-known result for zonoids (see \cite[Thm.~3.5.3]{Schneider} or \cite[Thm.~4.13]{Hug}).
In the following, measurability in a topological space $E$ always refers to the Borel $\sigma$-algebra on $E$. Let $\mu$ be a finite measure on $E$, let $E_0\subseteq E$  be  measurable and   $\mu(E\setminus E_0)=0$. In this case we say that $\mu$ is supported in $E_0$. If $g:E_0\to\R$ is a bounded and measurable function, then the integral of $g$ over $E$ with respect to $\mu$ is defined by choosing any measurable extension of $g$ to $E$ (and clearly this is independent of the particular extension).
The next lemma follows from the fact (applied with $S=E_0$) that if $S$ is a separable metric space, then probability measures with finite support on $S$ are weakly dense in the probability measures on $S$ (see the discussion on pages 72-73 of \cite{Billingsley}, Appendix III, Thm.~4 and the discussion after Thm.~5 on page 239 of the first edition (1968) of \cite{Billingsley} or \cite[Thm.~3]{Varadarajan1958}).

\begin{lemma}
\label{thm:discreteApprox}
  Let $(E,\rho)$ be a separable metric space and  $E_0\subseteq E$  a measurable subset.
  Let $\mu$ be a finite Borel measure on $E$ with $\mu(E\setminus E_0)=0$. Then there is a sequence of discrete Borel measures $\mu_j$, $j\in\N$, on $E$ with $\mu_j(E) = \mu(E)$ and $\mu_j(E\setminus E_0)=0$ such that if  $g \colon E_0\to\R$ is continuous and bounded, then
  \begin{equation}\label{eqneu1}
   \lim_{j\to\infty}\int g \diff \mu_j =\int g \diff \mu.
  \end{equation}
\end{lemma}

Let $\K_*$ be a Borel subset of $\K^n$. In the following, we always assume that $\K_*\neq\varnothing$. With the restriction of the Hausdorff metric, $\K_*$ is a separable metric space whose Borel $\sigma$-algebra coincides with the subspace $\sigma$-algebra induced on $\K_*$. In particular, we will be interested in the cases of the homothety invariant classes $\Pt^n$, the set of polytopes in $\R^n$, and the subclass $\Pt^n_k$   which is closed in $\K^n$.
For a Borel measure  $\nu$ on a separable metric space $E$, the support of $\nu$  is the complement of the largest open set on which $\nu$ vanishes and denoted by $\supp\nu$. Thus $\supp \nu$ is a closed set. If $\nu$ is a finite Borel measure on $\K_*$ with bounded support, then $\supp\nu$ is closed in $\K_* $ (but not compact in general). If $\K_*$ is closed in $\K^n$ and $\supp\nu$ is bounded, then $\supp\nu$ is compact.

If $\K_*$ is a closed and homothety invariant class of convex bodies (hence containing all singletons), then the Minkowski class $\mathfrak{M}(\K_*)$ consists of all finite Minkowski sums of convex bodies from $\K_*$ and all convex bodies in their closure.

Next we define the \emph{positive hull} of the support of a measure on $\K_*$.

\begin{definition}
{\rm
  Let $\mu$ be a probability measure on a Borel set $\varnothing\neq\K_*\subseteq\K^n$.  Then
  \[
    \mPos\mu \coloneqq \set*{ \sum_{i=1}^N \lambda_i L_i \;\colon N \in \N_0, \forall i \in [N]\colon \lambda_i \ge 0, L_i \in \supp\mu }
  \]
  denotes the set of nonnegative (finite) Minkowski combinations of convex bodies in $ \supp\mu$,  where $\supp\mu$ is defined with respect to the metric space $\K_*$. The empty sum is defined as $\{0\}$. If $\mu_1, \ldots, \mu_\ell$ are probability measures on $\K_*$, then
  \[
    \mPos\prn*{\mu_1, \ldots, \mu_\ell} \coloneqq \mPos\mu_1 \times \cdots \times \mPos\mu_\ell
  \]
  is the set of $\ell$-tuples with components in $\mPos\mu_1, \ldots, \mPos\mu_\ell$, respectively.}
\end{definition}

We provide a simple lemma. As usual, empty sums are interpreted as $0$ (or $\{0\}$ if sets in $\R^n$ are concerned). Recall that $\kappa_n$ is the volume of the unit ball $B^n$ and $\omega_n=n\kappa_n$ denotes its surface area. The mean width of a convex body $K\in\K^n$ can be expressed in the form
$$
w(K)=\frac{2}{\omega_n}\int_{\mathbb{S}^{n-1}} h_K(u)\, \mathcal{H}^{n-1}(\diff u)\ge 0
$$
with $w(K)>0$ if and only if $\dim K\ge 1$.

\begin{lemma}
\label{thm:diamSum}
  Let $\ell,n\in\N_0$ and  $A_1, \ldots, A_\ell\in\K^n$. Then
  \begin{equation}\label{eq:const}
    \sum_{i=1}^\ell \diam A_i \le \sqrt{\pi}n\, \diam \sum_{i=1}^\ell A_i
  .\end{equation}
\end{lemma}
\begin{proof} For the proof, we can focus on $n\ge 2$. Let $w(A)$ denote the mean width of $A\in\K^n$. Jung's inequality (or an obvious bound with $\sqrt{2}$ replaced by $2$) implies that $w(A)\le \sqrt{2}\diam A$. Moreover, since $A$ contains a segment of length $\diam A$, we have $2\diam A\le \omega_n\kappa_{n-1}^{-1} w(A)$.
Since the mean width is Minkowski additive, we get
\begin{align*}
    \sum_{i=1}^\ell \diam A_i&\le \frac{1}{2}\frac{\omega_n}{\kappa_{n-1}}\sum_{i=1}^\ell w(A_i)=\frac{1}{2}\frac{\omega_n}{\kappa_{n-1}}w\left(\sum_{i=1}^\ell A_i\right)
    \le \frac{\omega_n}{\kappa_{n-1}}\diam \sum_{i=1}^\ell A_i.
\end{align*}
The assertion follows since the Gamma function is increasing on $[1.5,\infty)$.
\end{proof}

 The representation \eqref{eq:refprokgen} in the following theorem  can be viewed as a specific version of Choquet's integral representation theorem (see \cite{Phelps} or \cite[Thm.~3.45 and Chap.~7]{LMNS10}), if combined with \cite[Thm.~3.4.2]{Schneider} (see also \cite[Thm. 4]{Berg69}). Thus it follows that the measure $\mu$  in Theorem \ref{thm:ptb-chargeneral} (b) can be chosen such that it is supported by the indecomposable convex bodies in $\K_*$.  We provide a direct argument for both directions of the following equivalence. The special case $\K_*=\Pt^n_k$ provides  a characterization of $k$-polyoids  and is stated as Corollary \ref{thm:ptb-char}. In the following, we always assume that $\varnothing\neq \K_*\subseteq \K^n$ is Borel measurable.

\begin{theorem}
\label{thm:ptb-chargeneral}
  Let $\varnothing\neq\K_*\subseteq\K^n$, $n\in\N_0$, be a homothety invariant closed class of convex bodies.   Then the following are equivalent.
  \begin{enumerate}[{\rm (a)}]
    \item $K\in\mathfrak{M}(\K_*)$.\label{it:pc1gen}
    \item There is a probability measure $\mu$ on $\K_*$ with bounded support such that
    \begin{equation}\label{eq:refprokgen} h_K = \int h_L \, \mu(\diff L).
    \end{equation}\label{it:pc3gen}
  \end{enumerate}
  If {\rm (\ref{it:pc3gen})} holds, then $K$ is the limit of a sequence in $\mPos \mu$.
\end{theorem}

\begin{proof}
The assertion is clear for $n=0$ or if $\diam K=0$. Hence, we can assume that  $n\ge 1$ and $\diam K>0$ in the following.

\medskip

\noindent
  \enquote{(\ref{it:pc1}) $\implies$ (\ref{it:pc3})}:
  Without loss of generality, $0 \in K\in\mathfrak{M}(\K_*)$. Hence  $K = \lim_{\ell\to\infty} Q_\ell$, where $Q_\ell = \sum_{i=1}^{m_\ell} Q_\ell^\tidx i$ with $Q_\ell^\tidx i\in\K_*$, $m_{\ell}>0$ and $\diam Q_\ell^\tidx i>0$. There are points $x_\ell \in Q_\ell$ and  $x_\ell^\tidx i \in Q_\ell^\tidx i$ with $x_\ell = \sum_{i=1}^{m_\ell} x_\ell^\tidx i$ such that $x_\ell \to 0$ as $\ell\to\infty$.
  Setting $P_\ell \coloneqq Q_\ell -  x_\ell$ and $P_\ell^\tidx i \coloneqq Q_\ell^\tidx i -  x_\ell^\tidx i\in \K_*$ for $\ell\in\N$ and $i \in m_\ell$, we have
  \[
  K = \lim_{\ell\to\infty} P_\ell = \lim_{\ell\to\infty} \sum_{i=1}^{m_\ell} P_\ell^\tidx i,\quad 0 \in P_\ell^\tidx i\in\K_*\quad\text{and}\quad \diam P_\ell^{\tidx i} > 0.
  \]

  The sequence $(\diam P_\ell)_\ell$ is bounded by some constant $d \in (0, \infty)$.
 For  $\ell\in\N$ and $i \in [m_\ell]$,  define  positive numbers
    \[
      d_\ell \coloneqq \diam P_\ell, \quad d_\ell^\tidx i \coloneqq \diam P_\ell^\tidx i, \quad e_\ell \coloneqq  \sum_{i=1}^{m_\ell} d_\ell^\tidx i, \quad
      c_\ell^\tidx i \coloneqq \frac{e_\ell}{d_\ell^\tidx i}
    ,\]
  and discrete probability measures $\mu_\ell$ on $\K_*$  by
  \[
    \mu_\ell \coloneqq \sum_{i=1}^{m_\ell} \frac{1}{c_\ell^\tidx i} \delta_{c_\ell^\tidx i P_\ell^\tidx i}, \quad\text{noting that } \sum_{i=1}^{m_\ell} \frac{1}{c_\ell^\tidx i} = 1
  .\]
  By construction and basic properties of support functions,
  \[
    h_{P_\ell} = \int h_P \,\mu_\ell(\diff P)
  .\]

  If $P \in \supp\mu_\ell$, then $P = c_\ell^\tidx i P_\ell^\tidx i$ for some $i \in [m_\ell]$, hence $0 \in P$ and, by Lemma \ref{thm:diamSum},
  \[
    \diam P = c_\ell^\tidx i d_\ell^\tidx i = e_\ell \le \sqrt{\pi}n\,d_\ell \le \sqrt{\pi}n\,d,
  \]
  so that $\supp\mu_\ell \subseteq S \coloneqq \set*{ L \in\K_*\;\colon L \subseteq \sqrt{\pi}nd\stdb }$. Since $\K_*$ is closed in $\K^n$, $S$ is compact.
  Thinking of the measures $\mu_\ell$ as measures on $S$ (with the restriction of the Hausdorff metric), a special case of Prokhorov's theorem  \cite[pp.~57--59]{Billingsley}  yields a subsequence $(\mu_{\ell_s})_s$ of $(\mu_\ell)_\ell$ that weakly converges to some probability measure $\mu$ on $\K_*$ which is also compactly supported in $S$.
  Hence, for all $u \in \R^n$,
  \[
    h_K(u) = \lim_{s\to\infty} h_{P_{\ell_s}}(u) = \lim_{s\to\infty} \int h_P(u) \,\mu_{\ell_s}(\diff P) = \int h_P(u)\, \mu(\diff P),\]
  since $P\mapsto h_P(u)$ is continuous and bounded on $S$. So $\mu$ has the desired property.

  \medskip

\noindent
  \enquote{(\ref{it:pc3}) $\implies$ (\ref{it:pc1})}:
  Let $\mu$ be a probability measure on $\K_*$ with bounded support such that \eqref{eq:refprokgen} holds. Let $E_0$ denote the support of
  $\mu$ with respect to the metric space $E=\K_*$.
  According to Lemma \ref{thm:discreteApprox},  $\mu$ is the weak limit of a sequence $\mu_\ell$ of discrete probability measures on $\K_*$ supported in $\supp\mu$. For all $\ell\in\N$, we define $K_\ell \in \K^n$ by
  \[
    h_{K_\ell} = \int h_P \, \mu_\ell(\diff P).
  \]
  By construction, $K_\ell \in \mPos\mu_\ell$ is a finite sum of convex bodies in $\K_*$.

  Since $\supp \mu$ is bounded, the function $P\mapsto h_P(u)$, $P\in E_0$, is bounded and continuous, for each $u\in \R^n$. Hence Lemma \ref{thm:discreteApprox}  ensures that, for each $u\in\R^n$,
  \[
    h_{K_\ell}(u) = \int h_P(u) \,\mu_\ell(\diff P) \to \int h_P(u) \,\mu(\diff P) = h_K(u) \quad (\ell\to\infty).
  \]
 This shows that $K_\ell \to K$ as $\ell\to\infty$ (with respect to the Hausdorff metric).
\end{proof}

\begin{corollary}
\label{thm:ptb-char}
  Let $K$ be a convex body in $\R^n$, $n\in\N_0$ and $k \in\N$. Then the following are equivalent.
  \begin{enumerate}[{\rm (a)}]
    \item $K$ is a $k$-polyoid.\label{it:pc1}
    \item There is a probability measure $\mu$ on $\Pt^n_k$ with compact support such that
    \begin{equation}\label{eq:refprok} h_K  = \int h_P  \, \mu(\diff P).
    \end{equation}\label{it:pc3}
  \end{enumerate}
  If {\rm (\ref{it:pc3})} holds, then $K$ is the limit of a sequence in $\mPos \mu$.
\end{corollary}

\begin{remark}
{\rm
  In view of Corollary \ref{thm:ptb-char}, a probability measure $\mu$ on $\Pt^n_k$ with compact support in $\Pt^n_k$  satisfying \eqref{eq:refprok} is called a \emph{generating measure} of the $k$-polyoid  $K$.
  Generating measures of polyoids are not uniquely determined (compare \cite[Rem.~3.2.15]{Schneider}). In the following, we will only use that for a given polyoid a generating measure exists.
}
\end{remark}

\begin{example} {\rm
    We describe the non-uniqueness by a simple example. Let  $e_1,e_2\in\R^2$ be the standard basis vectors.  Consider the intervals $I_1:=[0,e_1]$, $I_2:=[0,e_2]$ and $I_3:=[0,e_1+e_2]$. Let  $\conv$ denote the convex hull operator. Let
  $$
  P_1:=\conv(I_1\cup\{e_1+e_2\})\quad \text{and}\quad  P_2:=\conv\{e_1,e_1+e_2,2e_1+e_2\}.
  $$
  Then
  $$
  \mu_1:=\frac{1}{2}\left(\delta_{I_2+I_3}+ \delta_{I_1}\right),\quad  \mu_2:=\frac{1}{2}\left(\delta_{I_1+I_2}+\delta_{I_3}\right) \quad\text{and}\quad \mu_3:=\frac{1}{2}\left(\delta_{P_1}+ \delta_{P_2}\right)
  $$
  are three generating measures of the $4$-polyoid $P:=\frac{1}{2}(I_1+I_2+I_3)$, which in fact is also a zonoid (zonotope) with generating measure
  $$
  \mu_4:=\frac{1}{3} \left(\delta_{\frac{3}{2}I_1}+  \delta_{\frac{3}{2}I_2}+  \delta_{\frac{3}{2}I_3}\right).
  $$
  By adding to $P$ a suitable triangle, we get a $3$-polyoid which is not a zonoid and has two different generating measures.
  In the plane, examples of non-uniqueness can be easily constructed using Minkowski's existence theorem for polygons and the Minkowski additivity of the first area measure.}
\end{example}

Corollary \ref{thm:ptb-char} shows that polyoids can be characterized via the integral representation \eqref{eq:refprok} and as limits of sequences in the positive hull of a generating measure of the polyoid. The arguments in Section \ref{sec:3} are based on both types of description. In the following lemma, we show that a convex body the support function of which is given by a  more general integral representation is still the limit of a sequence of polytopes in the positive hull of a generating measure on $\Pt^n$. The lemma suggests the definition of a class of convex bodies that we will call macroids in Definition \ref{def:Macroid}. The argument for the implication \enquote{(b) $\implies$ (a)} of Lemma
\ref{thm:ptb-chargeneral} does not use any specific properties of the measurable subclass $\K_*\subseteq\K^n$. Therefore we  have the following lemma. Finally, we will choose $\K_*=\Pt^n$.

\begin{lemma}
\label{lem:ptb-con}
  Let $\varnothing\neq\K_*\subseteq\K^n$ be a  Borel set, $n\in\N_0$. Suppose that $\mu$ is a probability measure  on $\K_*$ with bounded support. Let $K\in\K^n$ be defined by
    \begin{equation}\label{eq:refprokb} h_K = \int h_P \, \mu(\diff P).
    \end{equation}
Then $K$ is the limit of a sequence in $\mPos \mu$.
\end{lemma}

\begin{definition}[Macroids]\label{def:Macroid}
    {\rm Let $\varnothing\neq\K_*\subseteq\K^n$ be a Borel set. A convex body $K$ in $\R^n$, $n\in\N_0$, for which there  is a probability measure   $\mu$  on $\K_*$ with bounded support such that \eqref{eq:refprokb} holds, is called a $\K_*$-\emph{macroid} with generating measure $\mu$. If $\K_*=\Pt^n$, we call $K$ a macroid with generating measure $\mu$.}
\end{definition}

\begin{remark}
{\rm Suppose that $K$ is a $\K_*$-macroid with generating measure $\mu$ on $\K_*$. We may extend $\mu$ trivially to all of $\K^n$. Then $\mu$ is a probability measure with bounded support (by definition) and (by Fubini's theorem)
$$
w(K)=\int w(Q)\, \mu(\diff Q)<\infty.
$$
The assumption that $\mu$ is a probability measure is not restrictive. To see this, note that if $\widetilde{\mu}$ is a Borel measure on $\K^n$ with $|\widetilde{\mu}|:=\widetilde{\mu}(\K^n)\in (0,\infty)$ and if $\widetilde{\mu}$ has bounded support, then
$$
\int h_Q\, \widetilde{\mu}(\diff Q)=\int h_Q\, \mu(\diff Q),
$$
where $\mu(\mathcal{A}):=|\widetilde{\mu}|^{-1}\widetilde{\mu}\left(|\widetilde{\mu}|^{-1}\mathcal{A}\right)$, for Borel sets $\mathcal{A}\subseteq \K^n$, defines a probability measure with bounded support.

For the present purpose, we could also replace the assumption of bounded support by an integrability assumption. To explain this statement, let
$\mu$ be a Borel probability measure on $\K^n$ such that
$$
0<\int w(Q)\, \mu(\diff Q)<\infty.
$$
The Steiner point of $K\in \K^n$ is defined by
$$
s(K):=\frac{1}{\kappa_n}\int_{\mathbb{S}^{n-1}} h_k(u)u\, \mathcal{H}^{n-1}(\diff u)
$$
and satisfies $s(K)\in \text{relint}\, K$ (see \cite[Sect.~1.7.1]{Schneider}).
Fubini's theorem yields
$$
s(K)=\int s(Q)\, \mu(\diff Q).
$$
Therefore we obtain
$$
h_{\frac{K-s(K)}{w(K)}}=\int h_P\, \mu^*(\diff P),
$$
where
$$
\mu^*(\mathcal{A}):=\frac{1}{w(K)}\int \mathbf{1}\left\{\frac{Q-s(Q)}{w(Q)}\in \mathcal{A}\right\}w(Q)\, \mu(\diff Q),
$$
for Borel sets $\mathcal{A}\subseteq \K^n$,
is a probability measure concentrated on
$$
\K^n_{0,1}:=\{L\in \K^n\colon w(L)=1,s(L)=0\}.
$$
In particular, $\mu^*$ has bounded support.
}
\end{remark}

\begin{remark}\label{rem:polmacrel}
    {\rm Each polyoid is a macroid, but not every macroid is a polyoid; for an example, see Appendix \ref{app:macroid-not-polyoid}. An explicit geometric characterization of the class of polyoids within the class of macroids remains to be discovered.}
    \end{remark}

    \begin{remark}{\rm
     An obvious motivation for introducing macroids is that Theorem \ref{corfin} is true in fact for the strictly larger class of macroids. An explicit example of a convex body that is not a macroid  is provided by a circular cone. This follows from   Proposition \ref{prop:indemacro}.
     }
     \end{remark}

     \begin{proposition}\label{prop:indemacro}
         Let $K \in \K^n$ be an indecomposable macroid. Then $K$ is a polytope.
     \end{proposition}

\begin{proof}
We may assume that $\dim K>0$ and
$$
h_K=\int h_Q\, \mu(\diff Q),
$$
where $\mu$ is a Borel probability measure on $\Pt^n$ with bounded support.
By Fubini's theorem, we have
$$
w(K)=\int w(Q)\, \mu(\diff Q),
$$
hence there is some $P\in \supp \mu$ with $w(P)>0$ (that is, $\dim P>0$) and $\mu(B(P,1/k))>0$ for all $k\in\N$, where $B(P,1/k)$ denotes a closed ball around $P$ with radius $1/k$ in $\Pt^n$ (or in $\K^n$) with respect to the Hausdorff metric $d$ on $\K^n$ (or its restriction to a subset).
    For  $k\in\N$, the convex body $K_k\in\K^n$ is defined by
    \[
       h_{K_k} \coloneqq \frac{1}{\mu(B(P, 1/k))} \int_{B(P, 1/k)} h_Q\, \mu(\diff Q)
    \]
    and satisfies $w(K_k)>0$.
    Then clearly $K_k\to P$ as $k\to\infty$ (with respect to the Hausdorff metric). Moreover, if $L_k\in\K^n$ is given by
        $$
        h_{L_k}:=\int_{B(P, 1/k)^\complement}h_Q\, \mu(\diff Q),
        $$
        then
    $$
    \mu(B(P,1/k))K_k+L_k=K.
        $$
    Since $K$ is indecomposable and $\dim K_k>0$, it follows that $K=c(k)K_k+x_k$, where
    $$
    c(k)=\frac{w(K)}{w(K_k)}\quad\text{and}\quad x_k\in\R^n.
    $$
    Since $K_k\to P$, we have $c(k)\to w(K)/w(P)>0$ and $x_k\to x_0\in\R^n$, as $k\to\infty$. Thus we arrive at $K=w(P)^{-1}w(K)P+x_0$, which shows that $K$ is a polytope.

     \end{proof}

     \begin{remark}{\rm
     Various types of mean section or projection bodies have been studied in integral and stochastic geometry. Starting from a convex body $K\subset\R^n$, the support function of a new mean body is defined as an integral average of the support functions of sections or projections of $K$, which is precisely the principle by which macroids are defined; see \cite{GW98,GW12,GW14,GHW17} and the literature cited there.

    Another special case of definition \eqref{eq:refprokb} is the convolution $\widetilde{\mu}\ast h_K$ of a probability (or finite) measure $\widetilde{\mu}$  on the rotation group ${\rm SO}_n$ and the support function of a fixed convex body $K\in\K^n$, as considered in \cite[Sects. 2 and 5]{GZ99}. In our notation, this reads
    \begin{align*}
    (\widetilde{\mu}\ast h_K)(u)&=\int h_K(\rho^{-1} u)\, \widetilde{\mu}(\diff \rho) =\int h_L(u)\, f_K(\widetilde{\mu})(\diff L),
    \end{align*}
    where $f_K(\widetilde{\mu})$ is the image measure of $\widetilde{\mu}$ under the map $f_K:{\rm SO}_n\to\K^n$,  $\rho\mapsto \rho K$.

    A general definition of a convex body as an integral average with respect to some measure on a suitable index set has been anticipated by Wolfgang Weil in \cite[(1)]{Weil76}, but then only the special case of zonoids has been explored in \cite{Weil76}.
    }
\end{remark}

\begin{remark}\label{rem:support-macroid}
    {\rm Let $u\in\mathbb{S}^{n-1}$. If $K$ is a macroid with generating measure $\mu$, then the support set  $F(K,u)$ of $K$  is a  macroid with generating measure $F_u(\mu)$, where $F_u:\Pt^n\to\Pt^n$, $P\mapsto F(P,u)$, is measurable and $F_u(\mu)$ is the image measure of $\mu$ under the  map $F_u$, that is,
    \begin{equation}\label{eqsupportset}
        h_{F(K,u)}=\int h_{F(P,u)}\, \mu(\diff P)=
    \int h_Q\, F_u(\mu)(\diff Q).
    \end{equation}
The measurability of $F_u$ follows from \cite[Thm.~12.2.6 (a) and Thm.~12.3.2]{SW}, since $F(K,u)=K\cap H(K,u)$, where $H(K,u)=u^\perp +h(K,u)u$ clearly depends continuously on $K$. Furthermore, note that $h_{F(K,u)}(x)=h_K'(u;x)$ by \cite[Thm.~1.7.2]{Schneider}, for $x\in\R^n$. Since
$t^{-1}|h_L(u+tx)-h_L(x)|\le R\|x\|$, for $t>0$ and $L\in \supp \mu\subseteq RB^n$ (and some $R>0$), the assertion follows from the dominated convergence theorem.
}
\end{remark}

\section{The characterization theorem}\label{sec:3}

We start by recalling various concepts of  criticality for finite sequences of subsets of $\R^n$. Recall that the cardinality of a finite set $I$ is denoted by  $|I|$. For a nonemtpy set $A\subseteq\R^n$, let $\vspan A$ denote the (linear) span of $A$ and $\pspan A=\vspan(A-A)$ the linear subspace parallel to the  smallest affine subspace containing $A$. By the dimension $\dim A\in\{0,\ldots,n\}$ of a set $A\neq\varnothing$ we mean the dimension of its affine span.

\begin{definition}\label{Def:critical}{\rm
    Let $\mathbfcal{A} = (A_1, \ldots, A_\ell)$, $\ell\in\N_0$, be a tuple of nonempty subsets of $\R^n$. Then $\mathbfcal{A}$ is called
    \begin{enumerate}[{\rm (i)}]
      \item \emph{supercritical} if $\dim \pspan \sum_{i\in I} A_i\ge \abs{I} + 2$ for all $\varnothing \ne I \subseteq \set*{1, \ldots, \ell}$.
    \item \emph{critical} if $\dim \pspan \sum_{i\in I} A_i\ge \abs{I} + 1$ for all $\varnothing \ne I \subseteq \set*{1, \ldots, \ell}$.
    \item \emph{semicritical} if $\dim \pspan \sum_{i\in I} A_i\ge \abs{I}$ for all $\varnothing \ne I \subseteq \set*{1, \ldots, \ell}$.
    \end{enumerate}
    }
\end{definition}

Note that here we deviate from the terminology used in \cite[Sect.~12]{SvH23+}, where a tuple of convex bodies satisfying (iii) in Definition  \ref{Def:critical} is called subcritical instead of semicritical. (Instead we reserve the notion of a subcritical tuple of sets for one that is not critical; see \cite{HugReichert23+}).

The various notions of criticality introduced above have useful properties some of which are discussed below. Each of the three notions is preserved by passing to a subtuple, taking permutations of the given tuple, replacing all sets by the same affine transformation or by individual translations, or if the sets are replaced by supersets. Supercriticality implies criticality, which in turn implies semicriticality. The empty tuple is supercritical. Moreover, if all sets  in an $\ell$-tuple $\mathbfcal{A}$ are full-dimensional, then $\mathbfcal{A}$ is supercritical if and only if $\ell \le n-2$ or $\ell=0$ (that is, $\mathbfcal{A}$ is the empty tuple).

\begin{lemma}
\label{thm:critSimple}
  Let $\ell \in \N$ and $\mathbfcal{A} = (A_1, \ldots, A_\ell)$ be a tuple of nonempty sets in $\R^n$.
  \begin{enumerate}[{\rm (a)}]
    \item Let $\mathbfcal{A}$ be critical and $A_{\ell + 1} \subseteq \R^n$ be nonempty. Then $(A_1, \ldots, A_{\ell + 1})$ is semicritical if and only if $\dim\pspan A_{\ell + 1}\ge 1$.\label{it:critSimple1}
    \item Let $\mathbfcal{A}$ be supercritical and $A_{\ell + 1} \subseteq \R^n$ be nonempty. Then $(A_1, \ldots, A_{\ell+1})$ is critical if and only if $\dim\pspan A_{\ell + 1}\ge 2$.\label{it:critSimple2}
  \end{enumerate}
\end{lemma}
\begin{proof}
(a)  Suppose that $A_{\ell + 1}$ has dimension at least one. Let $I\subseteq [n+1]$ be nonempty. We distinguish three cases.

  If $I\subseteq [\ell]$, then $\dim\pspan \sum_{i\in I}A_i\ge |I|+1\ge |I|$, since $\mathbfcal{A}$ is critical.

   If $I=\{\ell +1\}$, then  $\dim\pspan \sum_{i\in I}A_i\ge 1=|I|$, since $\dim\pspan A_{\ell + 1}\ge 1$.

   If $I=J\cup \{\ell +1\}$ and $\varnothing \neq J\subseteq [\ell]$, then
   $$
   \dim\pspan \sum_{i\in I}A_i\ge\dim\pspan \sum_{i\in J}A_i\ge |J|+1=|I|,
   $$
   where we used again that $\mathbfcal{A}$ is critical.

   Clearly, if $(A_1, \ldots, A_{\ell + 1})$ is semicritical then $\dim\pspan A_{\ell + 1}\ge 1$.

   \medskip

   The proof of (b) is similar.
\end{proof}

The following lemma  connects semicriticality of an $n$-tuple of convex bodies to the positivity of the mixed volume of  these convex bodies (see \cite[Theorem 5.1.8]{Schneider}).

\begin{lemma}\label{lem:critvol}
Let $\collc=(K_1,\ldots,K_n)$ be an $n$-tuple of convex bodies in $\R^n$. Then $\collc$ is semicritical if and only if \ $\V(\collc)>0$.
\end{lemma}

As pointed out before, mixed area measures can be extended to differences of support functions. If $g_1,g_2$ are differences of support functions and $\collc$ is an $(n-3)$-tuple of convex bodies in $\R^n$ (if $n\ge 3$), then we set $S_{g_1,g_2,\collc}:=S(g_1,g_2,\collc,\cdot)$. A similar convention applies in case just one of the bodies is replaced by a difference of support functions.

 The statement and proof of the following lemma is suggested by a similar result concerning zonoids; see \cite[Theorem~14.9]{SvH23+}.

 In the following, $\K_*\subseteq\K^n$ always denotes a measurable class of convex bodies.

\begin{lemma}\label{thm:extremalDecomposition}
  Assume that $n \ge 3$. Let $\collc$ be an $(n-3)$-tuple of  convex bodies in $\R^n$, and let $K\in\K^n$ be a $\K_*$-macroid with generating measure $\mu$. If $(K, \collc)$ is supercritical and $f$ is a difference of support functions with $\Su_{f, K, \collc} = 0$, then $\Su_{f, P, \collc} = 0$ \ for all $P \in \mPos\mu$.
\end{lemma}
\begin{proof}

  Let $P \in\mPos\mu$. Then there are $\ell\in\N_0$, $\lambda_1, \ldots, \lambda_\ell\ge 0$ and $L_1, \ldots, L_\ell \in\supp\mu$ such that $P = \sum_{i=1}^\ell \lambda_i L_i$ and
  \[
    \Su_{f, P, \collc} = \sum_{i=1}^\ell \lambda_i \Su_{f, L_i, \collc}
  .\]
  Note that this holds trivially with $\Su_{f, P, \collc}=0$ if $\ell=0$.
  So it suffices to prove that $\Su_{f, L, \collc} = 0$ for all $L \in \supp\mu$.

  By Fubini's theorem and basic properties of mixed area measures (see \cite[Sect.~5.1]{Schneider} or \cite[Sect.~4.1]{Hug}), which remain true in the case where  differences of support functions are admitted in some of the arguments of the mixed volumes and the mixed area measures,
  \begin{align}
    \int \V(f, f, L, \collc) \,\mu(\diff L)
    &= \frac{1}{n} \int \int h_L(u) \, \Su_{f, f, \collc}(\diff u) \, \mu(\diff L)\nonumber\allowdisplaybreaks\\
    & = \frac{1}{n} \int \int h_L(u) \, \mu(\diff L) \, \Su_{f, f, \collc}(\diff u)\nonumber\allowdisplaybreaks\\
    &= \frac{1}{n} \int h_K(u) \, \Su_{f, f, \collc}(\diff u)\nonumber \allowdisplaybreaks\\
    &= \V(f, f, K, \collc)\nonumber\\
    &= \frac{1}{n} \int f\,  \diff \Su_{f, K, \collc} = 0\label{eq:zerorel}
  .\end{align}

  If $L \in \K^n$ is a singleton,
   then $\V(f, f, L, \collc) = \V(f, f, 0L, \collc) = 0$ by translation invariance and multilinearity of $\V$.

  If $L \in \K^n$ is not a singleton, then $\V(K, K, L, \collc)>0$.  In fact, first we get $\dim K \ge 1+2=3$, since $(K, \collc)$ is supercritical.  By Lemma \ref{thm:critSimple} (b) it follows that $(K,K,\collc)$ is critical, but then $(L,K,K,\collc)$ is semicritical by Lemma \ref{thm:critSimple}(b) and since $\dim L\ge 1$. Hence the assertion follows from Lemma \ref{lem:critvol}.

  Since   $\Su_{f, K, \collc} = 0$, it follows from the extension of  \eqref{eqex} to differences of support functions that $ \V(f, K, L, \collc)=0$.
  Hence, by the General Alexandrov--Fenchel Inequality (GAFI)  we   get
  \[
    0 = \V(f, K, L, \collc)^2 \ge \V(f, f, L, \collc) \cdot \V(K, K, L, \collc),
  \]
  which implies that $\V(f, f, L, \collc)\le 0$.

 Since $\V(f, f, L, \collc)$ is continuous in $L \in \K^n$,  it follows from \eqref{eq:zerorel} that $\V(f, f, L, \collc) = 0$ for all $L \in\supp \mu$.

  Now let $L \in \supp\mu$. If $L$  is a singleton, then $\Su_{f, L, \collc} = 0$ by translation invariance and multilinearity of $\Su$. If $L$  is not  a singleton, then again $\V(K, K, L, \collc)>0$.   Moreover, $\V(f, K, L, \collc) = 0$   and $\V(f, f, L, \collc) = 0$, as shown above. Therefore, $\Su_{f, L, \collc} = 0$ is implied by \cite[Lem.~3.12 (a)]{SvH23+}.
\end{proof}

Next we compare how the smallest affine subspace containing a given $k$-polyoid with generating measure $\mu$ is related to the smallest affine subspace of a polytope from the positive hull of the support of $\mu$, if both affine subspaces are translated to the origin $0$.

\begin{lemma}\label{thm:mposSpan}
Let $n \in \N_0$. Let $K\in\K^n$ be a $\K_*$-macroid with generating measure $\mu$, and let $Q \in \mPos\mu$. Then $\pspan Q \subseteq \pspan K$.
\end{lemma}
\begin{proof} For $n=0$ the assertion is clear.
  Let $u \in \prn*{\pspan K}^\perp$ (the linear sub\-space orthogonal to $\pspan K$). Then
  \[
    \int (h_P(u) + h_P(-u))\, \mu(\diff P) = h_K(u) + h_K(-u) = 0
  .\]
Since the map $P\mapsto h_P(u)+h_P(-u)$, $P\in \K_*$, is continuous and nonnegative,  we get  $h_P(u) + h_P(-u) = 0$ for all $P \in \supp\mu$. Because $Q$ is a (nonnegative) Minkowski combination of such $P$, it follows that $h_Q(u) + h_Q(-u) = 0$. Hence, $u \in \prn*{\pspan Q}^\perp$. So $\prn*{\pspan K}^\perp \subseteq \prn*{\pspan Q}^\perp$, proving the claim.
\end{proof}

The proof of the next auxiliary result is inspired by \cite[Thm.~14.9]{SvH23+}.

\begin{lemma}
\label{thm:linGlue}
 Let $n\ge 3$. Let $\collc = (C_1, \ldots, C_{n-2})$ be an $(n-2)$-tuple of $\K_*$-macroids in $\R^n$ with generating measures $\mu_1, \ldots, \mu_{n-2}$. Let $f$ be a difference of support functions. Assume that $f$ is linear on $\supp\Su(\stdb, \collq,\cdot)$ whenever   $\collq = (Q_1, \ldots, Q_{n-2}) \in \mPos(\mu_1, \ldots, \mu_{n-2})$ with $\pspan Q_i = \pspan C_i$ \ for $i \in [n-2]$.
  Then $f$ is also linear on $\supp\Su(\stdb, \collc,\cdot)$.
\end{lemma}
\begin{proof}
By Lemma \ref{lem:ptb-con}, for   $i \in [n-2]$, there exists a sequence  $\widetilde C_i^\tidx 1, \widetilde C_i^\tidx 2, \ldots$  of sums of convex bodies in $\supp \mu_i$ that converges to $C_i$. Being an element of $\mPos\mu_i$, $\widetilde C_i^\tidx j$ satisfies $\pspan \widetilde C_i^\tidx j \subseteq \pspan C_i$ by Lemma \ref{thm:mposSpan}. On the other hand, $\widetilde C_i^\tidx j \to C_i$ implies that the reverse inclusion holds for all $j$ greater than or equal to some $q \in \N$. For   $i \in [n-2]$, define
  \[
    C_i^\tidx j \coloneqq \widetilde C_i^\tidx{q+j} + \frac{1}{j^2}\sum_{j'=1}^{j-1} \widetilde C_i^\tidx{q + j'},\quad j\in\N
  .\]
  Because $(d(0, \widetilde C_i^\tidx j))_j$ is bounded by some $c_i \in (0, \infty)$,
  \[
    d(C_i^\tidx j, \widetilde C_i^\tidx{q + j}) \le \frac{(j-1)c_i}{j^2} \to 0\quad \text{as }j\to\infty
  \]
  and
  \[
 \lim_{j\to\infty}   C_i^\tidx j = \lim_{j\to\infty} \widetilde C_i^\tidx {q+j} = C_i
  .\]
  Moreover,
  \[
    \pspan C_i = \pspan \widetilde C_i^\tidx{q + j} \subseteq \pspan C_i^\tidx j \subseteq \pspan C_i
  .\]

  For all $j \in \N$, we have $\collc^\tidx j \coloneqq (C_1^\tidx j, \ldots, C_{n-2}^\tidx j)\in \mPos(\mu_1, \ldots, \mu_{n-2})$.
  By assumption and since $\pspan C_i^\tidx j = \pspan C_i$ for $i \in [n-2]$, there is some $x_j\in\R^n$ such that $f = \left<x_j, \cdot\right>$ on $\supp\Su(\stdb, \collc^\tidx j,\cdot)$. By definition of $C_i^\tidx j$ and the multilinearity of $\Su$, we obtain
  \[
    \supp\Su(\stdb, \collc^\tidx j,\cdot) = \bigcup_{j'_1, \ldots, j'_{n-2} = 1}^j \supp\Su(\stdb, \widetilde C_1^\tidx{q + j'_1}, \ldots, \widetilde C_{n-2}^\tidx{q + j'_{n-2}},\cdot)
  .\]
  In particular,
  \begin{align}\label{eq:lg1}
    \supp\Su(\stdb, \collc^\tidx j,\cdot) \subseteq \supp\Su(\stdb, \collc^\tidx{j+1},\cdot) \quad \text{for all $j \in \N$}
  .\end{align}
  Hence, there is $p\in\N$ such that
  \[
    E \coloneqq \vspan \supp\Su(\stdb, \collc^\tidx p,\cdot) = \vspan \supp\Su(\stdb, \collc^\tidx j,\cdot) \quad \text{for all $j \ge p$}
  .\]
  Then for all $j \ge p$, $\left<x_p, \cdot\right>$ and $\left<x_j, \cdot\right>$ must agree on $E$ because they agree on $\supp\Su(\stdb, \collc^\tidx p,\cdot)$ and we obtain for all $j \ge p$,
  \[
    f = \left<x_p, \cdot\right> \quad \text{on $\supp\Su(\stdb, \collc^\tidx j,\cdot)$}
  .\]
  Because $\Su(\stdb, \collc^\tidx j,\cdot)$ weakly converges to $\Su(\stdb, \collc,\cdot)$, as $j\to\infty$, it follows that
  \[
    f = \left<x_p, \cdot\right> \quad \text{on $\supp\Su(\stdb, \collc,\cdot) \subseteq \cl\bigcup_{j=p}^\infty \supp\Su(\stdb, \collc^\tidx j,\cdot)$}
  ,\]
  which proves the assertion of the lemma.
\end{proof}

\begin{remark}\label{Remunnecessary}{\rm
In the proof of Lemma \ref{thm:linGlue}, we implicitly showed that if $\collc = (C_1, \ldots, C_{n-2})$ is an $(n-2)$-tuple of $\K_*$-macroids in $\R^n$, $n\ge 3$, with generating measures $\mu_1, \ldots, \mu_{n-2}$,  then there exists a sequence of $(n-2)$-tuples
$\collq^{(j)} = (Q_1^{(j)}, \ldots, Q_{n-2}^{(j)}) \in \mPos(\mu_1, \ldots, \mu_{n-2})$, $j\in\N$, such that $\pspan Q_i^{(j)} = \pspan C_i$ for $i\in [n-2]$ and $\collq^{(j)}\to \collc$ as $j\to\infty$. }
\end{remark}

We can now prove our main result. A crucial tool for our argument is the important special case of polytopes, which was already treated by Shenfeld and van Handel \cite[Thm.~8.1]{SvH23+}.  Recall that a convex body is said to be smooth if for each boundary point there is a unique supporting hyperplane passing through it. In particular, smooth convex bodies are full-dimensional.

\begin{theorem}
  \makeatletter\def\@currentlabelname{Characterization Theorem}\makeatother\label{thm:supercritical}
  Let $n\ge 2$. Let  $\collc = (C_1, \ldots, C_{n-2})$ be a supercritical $(n-2)$-tuple of macroids or smooth convex bodies in $\R^n$. Let $f$ be a difference of support functions. Then $\SufC = 0$ if and only if $f$ is linear on $\supp\Su(\stdb,\collc,\cdot)$.
\end{theorem}
\begin{proof}
Let $n=2$. Then the assumption states that $\Su_f=0$. If $f=h_K-h_L$ for $K,L\in \K^2$, this implies that $S(K,\cdot)=S(L,\cdot)$, hence $K=L+x$ for some $x\in\R^2$. This shows that $f$ is linear. Finally, note that $\supp \Su(B^n,\cdot)=\mathbb{S}^1$.

In the following, we assume that $n\ge 3$.
It is sufficient to prove the theorem in the case where $\collc$ is a supercritical tuple of macroids. The extension with the possible inclusion of smooth convex bodies follows immediately by an application of \cite[Cor.~14.3]{SvH23+}. For this, note that if a smooth convex body  in $\collc$  (which necessarily has full dimension) is replaced by the Euclidean unit ball (which is a zonoid and hence a polyoid), neither the condition $\SufC = 0$ nor the set $\supp\Su(\stdb,\collc,\cdot)$ are changed. Moreover, also the supercriticality of $\collc$ is not affected by replacing a smooth body by the unit ball.

Hence, it is sufficient in the following to consider  a supercritical $(n-2)$-tuple $\collc$ of macroids in $\R^n$ with $n\ge 3$.
First, we assume that $\SufC = 0$.
  Let $\collq = (Q_1, \ldots, Q_{n-2}) \in \mPos(\mu_1, \ldots, \mu_{n-2})$ be such that $\pspan Q_i = \pspan C_i$ for $i\in [n-2]$. So for all $I \subseteq [n-2]$, the tuple $(\collc_I, \collq_{I^\complement})$ is supercritical, where $I^\complement:=[n-2]\setminus I$. Based on the hypothesis $\SufC = 0$, Lemma \ref{thm:extremalDecomposition}  allows us to sequentially replace $C_i$ by $Q_i$ and to finally obtain $\Su_{f, \collq} = 0$.
  Since $\collq$ is a supercritical tuple of polytopes, $f$ is linear on $\supp\Su(\stdb, \collq,\cdot)$ by \cite[Thm.~8.1]{SvH23+}.
  Now the claim follows from Lemma \ref{thm:linGlue}.

  \medskip

For the reverse direction, we assume that $f$ is linear on $\supp\Su(\stdb,\collc,\cdot)$. Let $K\in\K^n$ be an arbitrary convex body. Then \cite[Lem.~7.6.15]{Schneider} (compare also  \cite[Lem.~8.11]{SvH23+}) implies that
$$
n\V(f,K,\collc)=\int f(u)\, \Su(K,\collc,\diff u)=0.
$$
By the symmetry of mixed volumes, we obtain
$$
0=\int h_K(u)\, \Su(f,\collc,\diff u),
$$
which yields $\Su(f,\collc,\cdot)=0$, since differences of support functions are dense in $C(\mbS^{n-1})$ (see e.g. \cite[Lem.~1.7.8]{Schneider}).
\end{proof}

Finally, we obtain a characterization of the equality cases of (AFI) for supercritical tuples of macroids and smooth bodies.

\begin{theorem}
\label{corfinmacroids}
Let $K, L \subset\R^n$ be convex bodies,  and let $\collc = (C_1, \ldots, C_{n-2})$ be a supercritical $(n-2)$-tuple of macroids or  smooth convex bodies in $\R^n$.
\begin{enumerate}[{\rm (a)}]
\item If $\V(K,L,\collc)=0$, then {\rm (AFI)} holds with equality and $K,L$ are homothetic.
\item Let $\V(K,L,\collc)>0$. Then {\rm (AFI)} holds with equality if and only if there are $a>0$ and $x\in\R^n$ such that $h_K=h_{aL+x}$ on $
\supp \Su(\stdb,\collc,\cdot)$.
\end{enumerate}
\end{theorem}

\begin{proof}
(a) If $\V(K,L,\collc)=0$, then $\V(K,K,\collc)\V(L,L,\collc)=0$, and hence (AFI) holds with equality.
By symmetry, we can assume that $\V(K,K,\collc)=0$. Then also $\V(K,K+L,\collc)=0$.

 If $K$ is a singleton, then $K,L$ are homothetic.

 If $K$ has dimension at least $1$, then $\dim(K+L)\le 1$, since otherwise $\dim(K+L)\ge 2$, $\dim(K)\ge 1$ and the
 assumed supercriticality of $\collc$ imply that $\V(K,K+L,\collc)>0$, a contradiction. Thus we get $\dim(K+L)\le 1$,
 in particular, $\dim(K)=1$ and $L$ is contained in a segment parallel to $K$. Hence again $K,L$ are homothetic.

 \medskip

 \noindent (b) Suppose that $\V(K,L,\collc)>0$. By \cite[Thm.~7.4.2]{Schneider} or \cite[Lem.~2.5]{SvH23+}, (AFI) holds with equality
 if and only if there is some $a>0$ such that
 $$
 \Su(K,\collc,\cdot)=\Su(aL,\collc,\cdot),
 $$
that is,
\begin{equation}\label{eq:equiv1}
\Su_{f,\collc}=0\quad \text{with }f=h_K-ah_L.
\end{equation}
Theorem \ref{thm:supercritical} implies that \eqref{eq:equiv1} holds if and only if here is some $x\in\R^n$ such that
\begin{equation}\label{eq:equiv3}
f=\langle x,\cdot\rangle\quad\text{on }\supp\Su(\stdb,\collc,\cdot),
\end{equation}
but clearly \eqref{eq:equiv3} is equivalent to
$$
h_K=h_{aL+x}\quad \text{on }\supp\Su(\stdb,\collc,\cdot),
$$
which proves the asserted equivalence.
\end{proof}

\appendix

\section{A macroid that is not a polyoid}\label{app:macroid-not-polyoid}

In this section, we construct an example of a macroid that is not a polyoid, thereby showing that the class of macroids is larger than the class of polyoids.

\subsection{Zonotope kernels of polytopes}

    Let $K, L, M \in \K^n$ be convex bodies. If $h_K - h_L = h_M$, then $K \ominus L \coloneqq M$ is called the \emph{Minkowski difference} of $K$ and $L$.

\begin{lemma}\label{thm:segment-summands}
  Let $K\in\K^n$ be a convex body and let $e, f$ be two linearly independent segments that are summands of $K$. Then $e + f$ is also a summand of $K$.
\end{lemma}
\begin{proof}
  To show that $e + f$ slides freely in $K$ (see \cite[Sect. 3.2]{Schneider} and in particular Theorem 3.2.2 there), it suffices to consider two-dimensional slices of $K$ parallel to $e + f$. Hence we can reduce to the case that $K$ is two-dimensional.

  Let $\pm u$ be the normals of $e$ and $\pm v$ the normals of $f$. As $F(e, \pm v)$ are trivial, $F(K \ominus e, \pm v)$ are translates of $F(K, \pm v)$. So translates of $f$ are not only contained in $F(K, \pm v)$ but also in $F(K \ominus e, \pm v)$. Then \cite[Thm.~3.2.11]{Schneider} yields that $f$ is a summand of $K \ominus e$. This completes the proof.
\end{proof}

\begin{lemma}\label{thm:zonotope-summand}
  The function
  $
    \zeta\colon \Pt^n \to \Pt^n
  $
  that maps a polytope to its unique largest (i.e.\ inclusion-maximal) zonotope summand, centrally symmetric around the origin, is well-defined. Every  zonotope summand of $P \in \Pt^n$ is a summand of $\zeta(P)$.
\end{lemma}
\begin{proof}
  We show that every polytope $P$ has a unique largest zonotope summand. Let $\mathcal{Z}(P)$ denote the nonempty set of origin centered zonotope summands of $P$.  First note that summands of polytopes are polytopes (see \cite[p.~157]{Schneider}) and polytopes that are zonoids are zonotopes (see \cite[Cor.~3.5.7]{Schneider}). Hence the set of all origin centered zonotopes that are summands of $P$ equals the set of all origin centered zonoids that are summands of $P$. The latter set is compact  as the intersection of a compact set (the set of centered zonoids having mean width less or equal the mean width of $P$) and a closed set (the set of summands of $P$).
   It follows that there is a $Z \in \mathcal{Z}(P)$ of maximum mean width. This $Z$ is inclusion-maximal in $\mathcal{Z}(P)$.

  Let $Y \in \mathcal{Z}(P)$.
  Then there are pairwise linearly independent $x_1, \ldots, x_k \in \stdsph$ and scalars $y_1, \ldots, y_k, z_1, \ldots, z_k \ge 0$ such that
  \[
    Y = \sum_{i=1}^k y_i [-x_i, x_i], \quad Z = \sum_{i=1}^k  z_i [-x_i, x_i]
  .\]
  Assume for a contradiction that $Y$ is not a summand of $Z$. Up to reordering of the indices, it follows that $y_1 > z_1$. Then $y_1 [-x_1, x_1]$ is a summand of $P$, but, as $Z$ is maximal, not a summand of
  \[
    P \ominus \sum_{i=2}^k z_i [-x_i, x_i].
  \]
  Let $\ell\in[k]$ be the largest index such that $y_1 [-x_1, x_1]$ is a summand of $\widetilde P \coloneqq P \ominus \sum_{i=2}^\ell z_i [-x_i, x_i]$. Then $l < k$ and $z_{\ell+1} [-x_{\ell+1}, x_{\ell+1}]$ is also a summand of $\widetilde P$. Now Lemma \ref{thm:segment-summands} shows that $y_1 [-x_1, x_1] + z_{\ell+1} [-x_{\ell+1}, x_{\ell+1}]$ is a summand of $\widetilde P$, but this contradicts the maximality of $\ell$. Hence, every $Y \in \mathcal{Z}(P)$ is a summand of $Z$, and there is only one maximal zonotope summand of $P$.
\end{proof}

Next, we aim to prove that $\zeta$ is measurable. We write $h(K,u)=h_K(u)$ for the support function of $K\in\K^n$ evaluated at $u\in\mathbb{S}^{n-1}$. We write $B(K,r)$ for a ball with center $K$ and radius $r\ge 0$  with respect to the Hausdorff metric $d$ on the space $\K^n$ of convex bodies.

\begin{lemma}
\label{thm:lowerSemicontinuous}
    Let $X$ be a separable metric space and $f \colon X \to \K^n$ a function such that for any $u \in \stdsph$ and $\lambda \in \R$,
    \[
        S_f(u, \lambda) \coloneqq \set*{ x \in X \colon h(f(x), u) \ge \lambda }
    \]
    is closed. Then $f$ is measurable.
\end{lemma}
\begin{proof}
    Fix some countable and dense set $Q \subseteq \stdsph$. Let $K \in \K^n$ and $r > 0$. By continuity of $h_L$ for every $L \in \K^n$,
    \[
    B(K, r) = \bigcap_{u \in Q } h(\cdot, u)^{-1}([h_K(u) - r, h_K(u) + r]).
    \]
    Taking the preimage under $f$, we get
    \[
    f^{-1}(B(K, r)) = \bigcap_{u \in Q}
    h(f(\cdot), u)^{-1}([h_K(u) - r, h_K(u) + r]).
    \]
    By hypothesis, $h(f(\cdot), u)$ is  a Borel set for every $u \in \stdsph$. Since $Q$ is countable, $f^{-1}(B(K, r))$
    is a Borel set as well.
    Because balls like $B(K, r)$ generate the Borel $\sigma$-algebra of $\K^n$, this shows that $f$ is measurable.
\end{proof}

\begin{lemma}\label{thm:zeta-measurable}
    $\zeta$ is a measurable function.
\end{lemma}
\begin{proof}
  We apply Lemma \ref{thm:lowerSemicontinuous}.
  Let $u \in \stdsph$ and $\lambda \in \R$. It suffices to show that
  \[
  S_\zeta(u, \lambda) = \set*{ P \in \Pt^n \colon h(\zeta(P), u) \ge \lambda }
  \]
  is closed. Let $(P_i)$ be a sequence in $S_\zeta(u, \lambda)$ that converges to $P \in \Pt^n$. Applying the Blaschke selection theorem to the bounded sequence $(\zeta(P_i) ) $, we find a subsequence $(Q_i)$ such that the sequence $(\zeta(Q_i))$ converges to a centered  zonoid $Z$ that is also a summand of $P$. Because summands of polytopes are polytopes  and polytopes that are zonoids are zonotopes, it follows that $Z$ is a zonotope. So $Z \in \zeta(P)$ and, in particular,
  \[
  h(\zeta(P), u) \ge h(Z, u) = \lim_{i\to\infty} h(\zeta(Q_i), u) \ge \lambda.
  \]
  So $P \in S_\zeta(u, \lambda)$, proving that the latter set closed. An application of Lemma \ref{thm:lowerSemicontinuous} concludes the proof.
\end{proof}

\subsection{Admissible sequences of polytopes}

    Let $K \subseteq \R^3$ be a convex body. A support set $F(K, u)$ will be called \emph{a singleton} or \emph{trivial} if it is zero-dimensional, \emph{an edge} if it is one-dimensional, and \emph{a facet} if it is two-dimensional. It should be observed that unless $K$ is a polytope, the current definition does not imply that the normal cone of $K$ at a  point in the relative interior of an edge is two-dimensional.

\begin{definition}\label{DefA3} \rm
    Let $(P_i)$ be a bounded sequence of (indecomposable) polytopes in $\R^3$ with the following properties:
    \begin{itemize}
        \item All facets are triangles.
        \item For every $i\in\N$, $P_i$ is $3$-dimensional.
        \item For every $i\in\N$, no two edges of $P_i$ have the same direction.
        \item If $i, j\in\N$ are distinct and $u$ is a facet normal of $P_i$, then $F(P_j, u)$ is trivial.
        \item If $\ell, i, j\in\N$ are distinct and $e, f, g$ are edges of $P_\ell, P_i, P_j$, then $e + f + g$ is $3$-dimensional. In particular, $e + f$ is $2$-dimensional.
        \item $K \coloneqq \sum_{i=1}^\infty P_i$ is a well-defined convex body.
    \end{itemize}
    We call such a sequence \emph{admissible} and $K$ its \emph{associated body}.
\end{definition}

\begin{remark}\label{rem:A1}\rm
Let $K_i$, $i\in\N$, and $K$ be convex bodies in $\mathcal{K}^n$. Then $K=\sum_{i=1}^\infty K_i$ holds (where the convergence of the partial sums is meant with respect to the Hausdorff metric) if and only if $h_K=\sum_{i=1}^\infty h_{K_i}$ (where the convergence holds pointwise, but then also uniformly on the unit sphere).
\end{remark}

\begin{remark}\label{rem:A2}\rm
Let $P_i,K\in\mathcal{K}^3$, $i\in\N$, be given as in Definition \ref{DefA3}.
Then $K$ has at most countably many extreme points. Items three, four and five imply that if $F(K,u)$ is an edge  of $K$,
then there is a unique $i\in\N$ such that $F(P_i,u)$ is an edge of $P_i$. In this situation, $F(K,u)$ is a translate of $F(P_i,u)$ and no other edge of any of the polytopes $P_j$, $j\neq i$, is parallel to $F(K,u)$. From  item four we conclude that if $F(K,u)$ is a triangular facet, then there is a unique $i$ such that $F(K,u)$ is a translate of $F(P_i,u)$. See Lemma \ref{thm:edges-of-faces} for further discussion.
\end{remark}

Recall that every summand of a polytope is a polytope (see \cite[p.~157]{Schneider}). For a polytope $Q\in\Pt^n$, we consider the convex cone
$$
\mathcal{S}(Q):=\set*{ P \in \Pt^n \given \exists R \in \Pt^n, \alpha > 0\colon Q = \alpha P + R }.
$$
The elements of $\mathcal{S}(Q)$ are called \emph{scaled summands of $Q$}.

\begin{lemma}\label{thm:polytope-measure-support}
    Let $Q \in \Pt^n$ be a polytope with macroid-generating measure $\mu$ on $\Pt^n$, that is,
    $$
    h_Q=\int h_P\, \mu(\diff P).
    $$
    Then $\supp\mu \subseteq \mathcal{S}(Q)$.
\end{lemma}
\begin{proof}~
    \textit{I.}
    Let $\beta > 0$ be a lower bound on the lengths of the edges of $Q$, and let $P\in \mathcal{S}(Q)$ be nontrivial. Then $\frac{\beta}{\diam P} P$ is a summand of $Q$, as we show first.

     Let $F(P, u)$ be an edge. Since $F(P, u)$ is a scaled summand of $F(Q, u)$, the latter must have an edge $e$ (which is also an edge of $Q$) homothetic to $F(P, u)$. The length of $F(P, u)$ is at most $\diam P$ and the length of $e$ is at least $\beta$, so $e$ contains a translate of $\frac{\beta}{\diam P} F(P, u)$.  Hence \cite[Thm.~3.2.11]{Schneider} implies that $\frac{\beta}{\diam P} P$ is a summand of $Q$.

    \textit{II.}
    The  set $\mathcal{S}(Q)$  is closed in $\Pt^n$, as we show next.

    Let $(P_i)$ be a sequence in $\mathcal{S}(Q)$ converging to some $P \in \Pt^n$. If $P$ is trivial, then $P \in \mathcal{S}(Q)$. Otherwise, there are a sequence of nontrivial polytopes $(P_i)$ and a sequence of polytopes $(R_i)$ such that $Q = \frac{\beta}{\diam P_i} P_i + R_i$, and $(R_i)$ must also converge to some $R \in \K^n$ such that $Q = \frac{\beta}{\diam P} P + R$. As $R$ is a summand of $P$, it must be a polytope. So $P \in \mathcal{S}(Q)$.

    \textit{III.}
    Assume for a contradiction that there is some $L \in \supp\mu\setminus \mathcal{S}(Q)$. Then ${\sf d} \coloneqq d(L, \mathcal{S}(Q)) > 0$ and $\lambda \coloneqq \mu(B(L, {\sf d}/2)) > 0$.
    Define convex bodies $L'$ and $R$ by
    $$h_{L'} \coloneqq \lambda^{-1} \int_{B(L, {\sf d}/2)} h_P \,\mu(\diff P)\quad\text{and}\quad h_R \coloneqq \int_{B(L, {\sf d}/2)^\complement} h_P \,\mu(\diff P)
    $$
    so that $Q = \lambda L' + R$. It follows that $R$ is a polytope and $L' \in \mathcal{S}(Q)$, and hence
    \[
      {\sf d} = d(L, \mathcal{S}(Q)) \le d(L, L') \le {\sf d}/2 < {\sf d},
    \]
    which is a contradiction.
\end{proof}

We write $\pspan A$ for the linear subspace parallel to the smallest affine subspace containing a given nonempty set $A\subseteq\R^3$.

\begin{lemma}\label{thm:exists-rational}
    For every edge $e$ of a polytope $P$ in $\R^3$, there are a normal $v \in \sphere{2}$ of $e$ and $u \in \mathbb{Q}^3\setminus\pspan \{e\}$ such that $u \perp v$ and $e = F(P, v)$.
\end{lemma}
\begin{proof}
    Let $U$ be the relatively open normal cone of the edge $e$. Then $\vspan U$ is two-dimensional, and $U$ is open in $\vspan U$. Choose $w\in \sphere{2} \cap U^\perp$, so that $\vspan U = w^\perp$.
    Let $\times$ denote the cross product. The continuous and surjective map
    $\R^3 \to \vspan U$, $x \mapsto w\times x$,
    maps the dense set $S \coloneqq \mathbb{Q}^3 \setminus \pspan \{e\} \subseteq \R^3$ to the dense set $\{w\} \times S \subseteq \vspan U$. Because $U \elminus0 \subseteq \vspan U$ is nonempty and open, there must be some $\tilde v \in (U\elminus0) \cap (\{w\} \times S)$. By construction, there is $u \in S$ such that $\tilde v = w \times u$.

    Now, $v \coloneqq \norm{\tilde v}^{-1}{\tilde v} \in \sphere{2}$ is an inner normal of $e$, i.e.\ $e = F(P, v)$, and is orthogonal to $u \in S = \mathbb{Q}^3 \setminus \pspan \{e\}$.
\end{proof}

In the following, we write $\pi_uK$ for the orthogonal projection of $K$ to $u^\perp$ for a vector $u \in \R^3 \elminus0$. Moreover, $\Su_1'(L,\cdot)$ denotes the first area measure of a convex body $L\subset u^\perp$ with respect to $u^\perp$ as the  ambient space.

\begin{lemma}\label{thm:edges}
    Let $(P_i)$ be an admissible sequence and $K$ its associated body together with a macroid-generating measure $\mu$.

    There is a set $\mathcal{M} \subseteq \Pt^3$ of full $\mu$-measure that satisfies the following property: If some $P \in \mathcal{M}$ has an edge with direction $v\in \sphere{2}$, then one of the $P_i$ also has an edge in direction $v$.
\end{lemma}
\begin{proof}
    We intend to use Lemma \ref{thm:exists-rational}. If $u \in \R^3 \elminus0$, then
    $$
\pi_{u^\perp} K=\sum_{i=1}^\infty \pi_{u^\perp} P_i,
    $$
    and hence the weak continuity and the Minkowski linearity of the area measures imply that
    \[
    \Su_1' (\pi_{u^\perp} K,\cdot)=\sum_{i=1}^\infty \Su_1'(\pi_{u^\perp} P_i,\cdot) .
    \]
    So $\Su_1'(\pi_{u^\perp} K,\cdot)$ is a discrete Borel measure (i.e., has countable support) in $u^\perp \cap \sphere{2}$. Denoting by
    \[
    \omega_u \coloneqq \set*{ v \in u^\perp \cap \sphere{2}  \colon \Su_1'(\pi_{u^\perp} K, \set*{v}) > 0 }
    \]
    the set of its atoms, we obtain from special cases of \cite[Thm.~2.23, Lem.~3.4]{HugReichert23+}  that
    \[
    0 = \Su_1'(\pi_{u^\perp} K, \omega_u^\complement) = \int \Su_1'(\pi_{u^\perp} P, \omega_u^\complement) \,\mu(\diff P).
    \]
    So the set
    \[
    \mathcal{M}_1 \coloneqq \bigcap_{u \in \mathbb{Q}^3\elminus 0} \set*{ P \in \Pt^3 \colon \Su_1'(\pi_{u^\perp} P, \omega_u^\complement) = 0 }
    \]
    has full $\mu$-measure.

    Since for each $u \in \mathbb{Q}^3\elminus0$ the set $\omega_u$ is countable, the set of pairs
    \[
    C \coloneqq \set*{ (u, v) \in (\mathbb{Q}^3\elminus0) \times \sphere{2} \colon  v \in \omega_u }
    \]
    is countable.
   Using the notation from Remark \ref{rem:support-macroid}, we define
    \[
    \mathcal{M}_2 \coloneqq \bigcap_{(u, v) \in C} F(\cdot, v)^{-1}(\supp F_v(\mu)).
    \]
    The set $\mathcal{M}_2$ has full $\mu$-measure. To see this, it is sufficient to consider $v\in\mathbb{S}^{n-1}$ and $P\in\Pt^3$ with $F(P,v)\notin \supp F_v(\mu)$ and to show that $P\notin\supp\mu$. By assumption, there is a neighbourhood $U'$ of $F(P,v)$ with $F_v(\mu)(U')=0$, hence $\mu(\{F (P',v):P'\in U'\})=0$. Since $\{F (P',v):P'\in U'\}$ is a neighbourhood of $P$, it follows that $P\notin\supp \mu$.

   Furthermore, for all $P \in \mathcal{M}_2$ and $(u, v) \in C$, Lemma \ref{thm:polytope-measure-support} shows that $F(P, v)$ is a scaled summand of $F(K, v)$.

    Now let $\mathcal{M} \coloneqq \mathcal{M}_1 \cap \mathcal{M}_2$. Assume $P \in \mathcal{M}$ and that $e$ is an edge of $P$. By Lemma \ref{thm:exists-rational}, there are $u \in \mathbb{Q}^3\setminus\pspan e$ and $v \in u^\perp \cap \sphere{2}$ such that $e = F(P, v)$. In particular, $F(\pi_{u^\perp} P, v)$ is nontrivial and $\Su_1'(\pi_{u^\perp} P, \set*{v}) > 0$. Since $P \in \mathcal{M}_1$, it follows that $v\in\omega_u$ and so $(u, v) \in C$. Since $P \in \mathcal{M}_2$, this implies that $e = F(P, v)$ is a scaled summand of the nontrivial support set $F(K, v)$, which is then either an edge or a parallelogram. If it is an edge, it has the same direction as $e$ and also is an edge of one of the $P_i$ and we are done. If it is a parallelogram, then one of the sides of the parallelogram must have the same orientation as $e$, and one of the $P_i$ has an edge with this orientation. This concludes the proof.
\end{proof}

\begin{lemma}\label{thm:support-set-summands}
  Let $(P_i)$ be an admissible sequence and $K$ its associated body together with a macroid-generating measure $\mu$. Then there is a set $\mathcal{M}'$ of full $\mu$-measure such that for each $P \in \mathcal{M}'$ and for each $u \in \sphere{2}$, $F(P, u)$ is a scaled summand of $F(K, u)$.
\end{lemma}
\begin{proof}~
    (I)
    Let $P \in \mathcal{M}$, where $\mathcal{M}$ is as in the statement of Lemma \ref{thm:edges}.
    If $F(P, u)$ is trivial, so is the claim. If $F(P, u)$ is an edge, then by Lemma \ref{thm:edges}, there is $i \in \N$ such that $F(P, u)$ is homothetic to the edge $F(P_i, u)$ and hence a scaled summand of $F(K, u)$.

   (II)
    Now consider the case that $P\in \mathcal{M}$ and $F(P, u)$ is a facet. The edges of the polytopes $P_i$, $i \in \N$, together have only countably many directions. Denote the countable set of these directions by $A \subset \sphere{2}$. The facet $F(P, u)$ is incident to (at least) two edges with linearly independent directions $v, w \in A$ that determine the facet normal $u$ up to sign $\sigma \in \set*{-1, 1}$ via
    \[
        \phi\colon \set*{-1, 1} \times \set*{ (a, b) \in A^2 \given a \ne \pm b }, \quad (\sigma, u, v) \mapsto \sigma \frac{u \times v}{\norm{u \times v}}.
    \]
    So the facet normals of $P$ are contained in the countable image of $\phi$, which is independent of the choice of $P$.

    For each $u \in \sphere{2}$, the set $F(K, u)$ is a polytope. Consider the set of full $\mu$-measure
    \[
      \mathcal{M}_3 \coloneqq \bigcap_{u\in\im\phi} F(\cdot, u)^{-1}(\supp F_u(\mu)).
    \]
    If $P$ is also in $\mathcal{M}_3$, then by Lemma \ref{thm:polytope-measure-support} and $u \in \im\phi$, the support set $F(P, u)$ is a scaled summand of $F(K, u)$.

     Now the assertion follows from (I) and (II) with $ \mathcal{M}' \coloneqq \mathcal{M} \cap \mathcal{M}_3$.
\end{proof}

\subsection{Unique decomposability}

\begin{lemma}\label{thm:edges-of-faces}
    Let $(P_i)$ be an admissible sequence and $K$ its associated body together with a macroid-generating measure $\mu$. Let $e$ be an edge of $F(K, u)$ for some $u \in \sphere{2}$. Then there is a unique $i\in\N$ such that $F(P_i, u)$ has an edge homothetic to $e$, $e$ is in fact a translate of $F(P_i,u)$, and this edge is unique among the edges of $P_i$.
\end{lemma}
\begin{proof}
  The uniqueness statements immediately follow from the properties of admissible sequences  $(P_i)$. Note that an edge of $F(K,u)$ need not be an edge of $K$ as defined here.

  If $F(K, u)$ is a singleton, it does not have any edges.

  If $F(K, u)$ is an edge, then there is $i\in\N$ such that $F(P_i, u)$ is a translate of $F(K, u)$.

  If $F(K, u)$ is a triangle, then there is $i\in\N$ such that $F(P_i, u)$ is a translate of $F(K, u)$. So  a translate of $e$ is an edge of $F(P_i, u)$.

  If $F(K, u)$ is a parallelogram, then there are unique $i, j\in\N$ with $i\neq j$ such that $F(P_i,u)$ is a translate of an edge of $P_i$, $F(P_j,u)$ is a translate of an edge of $P_j$ and $F(P_i, u) + F(P_j, u)$ is a translate of $F(K, u)$. So  $e$
  is either a translate of $F(P_i, u)$ or of $F(P_j, u)$.
\end{proof}

\begin{lemma}\label{thm:no-zonotope-summands}
    Let $(P_n)$ be an admissible sequence and $K$ its associated body together with a macroid-generating measure $\mu$.

    Then the polytopes with a nontrivial zonotope summand are contained in a $\mu$-zero set $\mathcal{N}$.
\end{lemma}
\begin{proof}
    Recall the measurable function $\zeta$ from Lemmas \ref{thm:zonotope-summand} and \ref{thm:zeta-measurable}. The macroid $Z$ generated by $\zeta(\mu)$, the image measure of $\mu$ under the map $\zeta$, is a zonoid and summand of $K$, implying
    \[
        \Su_2(Z,\cdot) \ll \Su_2(K,\cdot) = \sum_{n=1}^\infty \sum_{m=1}^\infty \Su(P_n, P_m,\cdot).
    \]
    Because the right-hand side is a discrete measure, so is $\Su_2(Z,\cdot)$. If we can show that $Z$ has no facets, then $\Su_2(Z,\cdot)$ is a discrete measure without atoms, hence zero. Then $Z$ is at most one-dimensional. If we can also show that $Z$ has no edges, then $Z$ must be trivial,
    and the set $\mathcal{N}$ of polytopes $P$ with nontrivial $\zeta(P)$ is a $\mu$-zero set, proving the claim.

    It remains to show that $Z$ has no facets or edges. We aim at a contradiction and assume that $F(Z, u)$ is a facet or an edge.

    Since $F(Z, u)$ is a summand of the polytope $F(K, u)$, it has an edge $e$ that is homothetic to an edge of $F(K, u)$. Because $Z$ is centrally symmetric around the origin, $F(Z, -u) = -F(Z, u)$ also has an edge that is a translate of $e$, and therefore $F(K, -u)$ has an edge that is homothetic to $e$. Hence, $F(K, u)$ and $F(K, -u)$ both contain an edge homothetic to $e$. By Lemma \ref{thm:edges-of-faces}, and especially the uniqueness statement, there is $i\in\N$ such that $F(P_i, u)$ and $F(P_i, -u)$ intersect in the very same edge homothetic to $e$. But this contradicts $P_i$ being $3$-dimensional, and $Z$ cannot have edges or facets.
    This completes the proof.
\end{proof}

\begin{lemma}\label{thm:unique-decomposability}
    Let $(P_i)$ be an admissible sequence and $K$ its associated body together with a macroid-generating measure $\mu$.

    Then $\mu$ is supported in translates of $\operatorname{pos} ((P_i)_i)$, the set of translates of finite positive (Minkowski) combinations of polytopes from the sequence $(P_i)_{i\in\N}$.
\end{lemma}
\begin{proof}
  By Lemmas \ref{thm:support-set-summands} and \ref{thm:no-zonotope-summands},
    $\mu$ is supported in the polytopes $P$ that have no nontrivial zonotope summand and such that for all $u \in \sphere{2}$, the support set $F(P, u)$ is a scaled summand of $F(K, u)$. Let $\mathcal{M}_4$ be a set of such polytopes of full $\mu$-measure, and let $P \in \mathcal{M}_4$. For the proof, we may assume that $P$ is nontrivial.

\begin{enumerate}[(i)]
    \item All facets of $K$ are triangles and parallelograms. The only scaled summands of a triangle are homothets of that triangle; the only scaled summands of a parallelogram are (possibly degenerate) parallelograms with the same edge directions but possibly different proportions. So all facets of $P$ are of this kind.
    \item Let $u\in\sphere2$ be such that $F(P, u)$ is a triangular facet. Then $F(K, u)$ is homothetic to $F(P,u)$, that is, there are unique $\alpha_u > 0$ and $t_u \in \R^n$ such that $F(P, u) = \alpha_u F(K, u) + t_u$. Moreover, there are unique $i = i(u)\in\N$ and $t_u'\in\R^n$ such that $F(K, u)$ is a translate of $F(P_i, u)$ and
    $F(P, u) = \alpha_u F(P_i, u) + t_u'$.

    Also note that there are at most two triangular facets of $P$  that have an edge parallel to a fixed direction; otherwise, there would be an $i\in\N$ such that $P_i$ also had more than two such facets, contradicting the hypothesis that $P_i$ has at most one edge parallel to  a given direction.
     \item We observe that $\dim P=3$. Recall that $P$ is nontrivial. If $\dim P=1$, then $P$ is a segment, which is a zonotope, a contradiction. If $\dim P=2$, then $P$ is a triangle or a non-degenerate parallelogram, which is a zonotope. The latter is excluded. Hence $P$ is a triangle with $P=F(P,v)=F(P,-v)$ for a unit vector $v$. Let $e$ be an edge of $P$. Then $P_{i(v)}$ and $P_{i(-v)}$ both contain an edge parallel to $e$. Hence, $i \coloneqq i(v) = i(-v)$ and $F(P_i, v) = F(P_i, -v)$ and thus $\dim P_i=2$, a contradiction.
    \item Let $G$ be the graph with the edges of $P$ as $G$-vertices, where two edges, i.e.\ $G$-vertices, are connected if and only if they are opposite edges in a parallelogram facet of $P$. Since every edge is only part of two facets, the maximum degree of a $G$-vertex, i.e.\ an edge of $P$, is two. The connected components of $G$ are cycles or chains.

    Let us first make sure that no cycles can occur. Assume that the edge $e$ of $P$ with direction $u$ is part of a cycle. Then $\pi_{u^\perp} P$ is a convex polygon and $\pi_{u^\perp} e$ is one of its vertices. The two edges incident to $\pi_{u^\perp} e$ are projections of parallelograms that connect $e$ to the two neighbors of $e$ in $G$, and an induction shows that all support sets of $\pi_{u^\perp} P$ either are projections of edges parallel to $e$ or parallelograms connecting two such edges.  For the sake of applying \cite[Thm.~3.2.22]{Schneider}, let $F(e, v)$ be an edge of $e$. In this case, $v \in u^\perp$, and so $e$ is a summand of $F(P, v)$. Then \cite[Thm.~3.2.22]{Schneider} guarantees that $e$ is a summand of $P$, in contradiction to $P$ having no nontrivial zonotope summand. Therefore, the connected component of any edge $e$ of $P$ is a chain $e_1 - \cdots - e_k$ of $e$-translates.

    The endpoints of this chain must be edges of two triangular facets of $P$. By (ii), there can be no other chain with edges parallel to $e_1, \ldots, e_k$. So if $f$ is an edge parallel to $e$, then $f = e_j$ for some $j\in[k]$ and $f$ is a translate of $e$. Moreover, for any edge $e$ of $P$ there are exactly two triangular facets of $P$ with an edge parallel to $e$.
    \item Let $u, v\in\sphere2$, and $i \coloneqq i(u)$ as in (ii), such that $F(P, u)$ is a triangle and $F(P_i, v)$ is a facet adjacent to $F(P_i, u)$ via an edge $e$. By (iii), there is exactly one $w\in\sphere2$ besides $u$ such that $F(P, w)$ is a triangle with an edge parallel to $e$. By (ii), $F(P_i, w) \ne F(P_i, u)$ is then also a triangle with an edge parallel to $e$. Because $P_i$ contains no other edge parallel to $e$, it follows that $v = w$. So we have
    \[
    F(P, u) = \alpha_u F(P_i, u) + t_u, \quad F(P, v) = \alpha_v F(P_i, v) + t_v.
    \]
    By (iii), all edges of $P$ parallel to $e$ are translates of each other, hence it follows that   $\alpha_u = \alpha_v$.
    \item Let $u, v\in\sphere2$, and $i \coloneqq i(u)$ as in (ii), such that $F(P, u)$ and $F(P_i, v)$ are triangles. Then the triangles $F(P_i, u)$ and $F(P_i, v)$ are connected via a chain of neighboring facets. Iteration of (v) shows that $F(P, v)$ is a triangle and $\alpha_u = \alpha_v$. So $\alpha_u$ only depends on $P$ and $i(u)$, and we set $\alpha_{i(u)}(P) \coloneqq \alpha_u > 0$. If $i\in\N$ and $P$ contains no triangular facet $F(P, w)$ with $i(w) = i$, then we set $\alpha_i(P) \coloneqq 0$.

    For each $i\in\N$, there are uncountably many edge normals of $P_i$ but only countably many facet normals of $K$. Let $u\in\sphere2$ be such that $F(P_i, u)$ is an edge and $F(K, u)$ is not a facet and in fact a translate of $F(P_i, u)$. Then for each $P \in \mathcal{M}_4$, $F(P, u)$ is a summand of $F(K, u)$ and satisfies $\V(F(P, u)) = \alpha_i(P) \V(F(P_i, u))$. This shows that  $\mathcal{M}_4 \ni P \mapsto \alpha_i(P)$ is measurable,
    \[
    \V(F(K, u)) = \int_{\mathcal{M}_4} \V(P, u) \mu(\diff P) = \int_{\mathcal{M}_4} \alpha_i(P) \mu(\diff P) \V(F(P_i, u))
    \]
 and thus we get
 \begin{equation}\label{eq:A14.int}
 \int_{\mathcal{M}_4} \alpha_i(P)\,\mu(\diff P) = 1.
 \end{equation}
    \item Let $\widetilde P \coloneqq \sum_{i=1}^\infty \alpha_i(P)P_i$, involving only finitely many nonzero summands. Note that $\dim\widetilde P=3$, since $\alpha_i(P)>0$ for some $i\in\N$.  Every facet of $P$ or $\widetilde P$ is either triangular or a parallelogram.
The preceding items show that the triangular facets of $P$ are translates of the triangular facets of $\widetilde P$, and vice versa.
\item
It remains to consider the parallelogram facets.
Let $u\in\sphere2$. If $F(K, u)$ is not a parallelogram, it is a singleton, an edge or a triangle. For all $P\in\mathcal{M}_4$, neither $F(P, u)$ nor $F(\widetilde P, u)$ is then a parallelogram.
\item
We consider the situation from (viii). From now on, we assume that $F(K, u)$ is a parallelogram. We choose $v, w\in\sphere2\cap u^\perp$ such that the edges of $F(K, u)$ are $F(F(K, u), \pm v)$ and $F(F(K, u), \pm w)$. There are unique distinct $i, j \in \N$ such that $F(F(K, u), \pm v)$ are translates of $F(P_i, u)$ and $F(F(K, u), \pm w)$ are translates of $F(P_j, u)$. Let $P \in \mathcal{M}_4$. Then   $F(\widetilde P, u)$ is a translate of
\[
\alpha_i(P)F(P_i, u) + \alpha_j(P)F(P_j, u),
\]
which might be a singleton, an edge or a parallelogram.
On the other hand,  $F(P, u)$ is a translate of
\[
    F(F(P, u), v) + F(F(P, u), w).
\]
\item We consider the situation from (ix) and aim to show that $F(P, u)$ and $F(\widetilde{P},u)$ are translates of each other. In the current item, we show that the conclusion holds at least if $P$ is taken from  a subset of $\mathcal{M}_4$ of full measure. The argument will be completed in (xi).

Let $P\in\mathcal{M}_4$. If $F(F(P, u), v)$ is not a singleton, then it is an edge parallel to $F(P_i, u)$. Hence it must be a translate of $\alpha_i(P)F(P_i, u)$. In either case,
\begin{equation}\label{eq:A14eq2}
    \V(F(F(P, u), v)) \le \alpha_i(P) \V(F(P_i, u)).
\end{equation}
Relation \eqref{eq:A14eq2}  can be used to bound the integrand in
\[
    \V(F(P_i, u)) = \V(F(F(K, u), v)) = \int_{\mathcal{M}_4} \V(F(F(P, u), v)) \,\mu(\diff P).
\]
But then \eqref{eq:A14.int}  from (vi) implies that equality must hold in \eqref{eq:A14eq2}, for $\mu$-almost all polytopes $P\in\mathcal{M}_4$. Hence, $F(F(P, u), v)$ is a translate of $\alpha_i(P) F(P_i, u)$, and a similar argument shows that $F(F(P, u), w)$ is a translate of $\alpha_j(P) F(P_j, u)$, for $\mu$-almost all $P\in\mathcal{M}_4$. So there is a measurable set $\mathcal{M}_5(u) \subseteq \mathcal{M}_4$ of full $\mu$-measure such that for all $P\in\mathcal{M}_5(u)$, a translate of $F(P, u)$ is
\[
\alpha_i(P) F(P_i, u) + \alpha_j(P) F(P_j, u),
\]
and hence also a translate of $F(\widetilde P, u)$.
\item Finally, set $\mathcal{M}_5 \coloneqq \bigcap_u \mathcal{M}_5(u)$, where we take the countable intersection over all normals of parallelogram facet normals of $K$. Then $\mathcal{M}_5$ is a measurable set of full $\mu$-measure and for all $P \in \mathcal{M}_5$ and $u \in \sphere2$, $F(P, u)$ is a parallelogram if and only if $F(\widetilde P, u)$ is, and in this case both are translates of each other: When $F(K, u)$ is a parallelogram, it follows from (x) that $F(P, u)$ and $F(\tilde P, u)$ are translates, and if it is not, neither of them is a parallelogram due to (viii).
\end{enumerate}
The proof is concluded by an application of Minkowski's uniqueness theorem for area measures of convex polytopes.
\end{proof}

\begin{lemma}\label{thm:prime-summand}
    Let $(P_i)$ be an admissible sequence and $K$ its associated body together with a  macroid-generating measure $\mu$.

    Then for all $i\in\N$ there is $P\in\supp\mu$ such that $P_i$ is a scaled summand of $P$.
\end{lemma}
\begin{proof}
  Let $u$ be the normal of a (necessarily triangular) facet of $P_i$. Then $F(K, u)$ is a translate of this triangular facet. Clearly, $\mu$ is concentrated on  $\supp\mu$ and, according to Lemma \ref{thm:unique-decomposability}, on the set of translates of all finite positive  combinations of the $P_j$, $j\in\N$. By  \eqref{eqsupportset}  we have
  \[
    h_{F(K, u)} = \int h_{F(P, u)}\, \mu(\diff P) ,
  \]
  hence there is $P \in \operatorname{pos}((P_j)_j) \cap \supp\mu$ such that $F(P, u)$ is nontrivial. This can only be the case if $P_i$ is a scaled summand in the finite positive  combination defining $P$, as $F(P_j, u)$ is trivial for all $j \ne i$.
\end{proof}

\begin{theorem}
    Let $(P_i)$ be an admissible sequence and $K$ its associated body.
    If for all $i\in\N$ the body $P_i$ has at least $i$ vertices, then $K$ is not a polyoid, although it is a macroid.
\end{theorem}
\begin{proof}
  Assume that $K$ is a $k$-polyoid with a generating measure $\mu$ supported in the space of $k$-topes. Lemma \ref{thm:prime-summand} shows that $P_{k+1}$ is a scaled summand of some $P\in\supp\mu$. But then $P$ is not a $k$-tope (see, e.g., \cite[Lem.~2.3]{DP19}), in contradiction to the property $\supp\mu \subseteq \Pt^3_k$ of $\mu$. So $K$ is not a $k$-polyoid for any $k\in\N$.
\end{proof}

\begin{remark} \rm
    Let $(P_i)_{i\ge 4}$ be a bounded sequence of polytopes such that $P_i$ is a $3$-dimensional $i$-tope having only triangular faces with no edge direction occurring twice. When we apply independent uniform random rotations to each of the $P_i$, we obtain almost surely an admissible sequence. This way, we can construct a macroid that is not a polyoid.
\end{remark}

\bigskip

\noindent
\textbf{Acknowledgements.}
D. Hug was supported by DFG research grant HU 1874/5-1 (SPP 2265). The authors are grateful to Ramon van Handel for helpful comments on an earlier version of the manuscript.


\bigskip

\vspace{2cm}

\noindent
Authors' addresses:

\bigskip

\noindent
Daniel Hug, Karlsruhe Institute of Technology (KIT), Institute of Stochastics, D-76128 Karlsruhe, Germany, daniel.hug@kit.edu

\medskip

\noindent
Paul Reichert, Karlsruhe, Germany, math@paulr.de


\begin{thebibliography}{99}

\bibitem{AF1937}
Aleksandr D.~Alexandrov. A. D. Alexandrov selected works. Part I. In: 1st edition.
Vol. 4. Classics of Soviet Mathematics. London: Yu. G. Reshetnyak
and S. Kutateladze, 1996. Chap. To the theory of mixed volumes of convex
bodies. Part II: New inequalities for mixed volumes and their applications.

\bibitem{Berg69}
Christian Berg. Shephard's approximation theorem for convex bodies and the Milman theorem. Math. Scand. 25 (1969), 19--24.


\bibitem{Billingsley}
Patrick Billingsley. Convergence of Probability Measures. Second edition. Wiley
Series in Probability and Statistics: Probability and Statistics. A Wiley-Interscience Publication.
John Wiley \& Sons, Inc., New York, 1999, pp. x+277.


\bibitem{DP19}
Antoine Deza and Lionel Pournin. Diameter, decomposability, and Min\-kow\-ski sums of polytopes.
Canad. Math. Bull. 62 (2019), no. 4, 741--755.

\bibitem{EKMS19}
Dario Cordero-Erausquin, Bo'az Klartag,  Quentin Merigot and Filippo Santambrogio. One more proof of the Alexandrov-Fenchel inequality. C. R. Math. Acad. Sci. Paris 357 (2019), no. 8, 676--680.

\bibitem{GW98}
Paul Goodey, Markus Kiderlen and Wolfgang Weil. Section and projection means of convex bodies. Monatsh. Math. 126 (1998), no. 1, 37--54.

\bibitem{GW12}
Paul Goodey and  Wolfgang Weil. A uniqueness result for mean section bodies. Adv. Math. 229 (2012), no. 1, 596--601.

\bibitem{GW14}
Paul Goodey and  Wolfgang Weil. Sums of sections, surface area measures, and the general Minkowski problem. J. Differential Geom. 97 (2014), no. 3, 477--514.

\bibitem{GHW17}
Paul Goodey, Daniel Hug and  Wolfgang Weil. Kinematic formulas for area measures. Indiana Univ. Math. J. 66 (2017), no. 3, 997--1018.


\bibitem{GZ99}
Eric Grinberg and Gaoyong Zhang.
Convolutions, transforms, and convex bodies.
Proc. London Math. Soc. (3) 78 (1999), no. 1, 77--115.




\bibitem{HugReichert23+}
Daniel Hug and Paul A. Reichert. The support of mixed area measures involving a new class of convex bodies. ArXiv.

\bibitem{Hug}
Daniel Hug and Wolfgang Weil. Lectures on Convex Geometry. Vol. 286.
Graduate Texts in Mathematics. Springer, Cham, 2020, pp. xviii+287.

\bibitem{LMNS10}
 Jaroslav Luke\v{s}, Jan Mal\'{y},  Ivan Netuka and Ji\v{r}\'{i} Spurn\'{y}.  Integral representation theory. Applications to convexity, Banach spaces and potential theory. De Gruyter Studies in Mathematics, 35. Walter de Gruyter \& Co., Berlin, 2010. xvi+715 pp.

\bibitem{Phelps}
 Robert R.~Phelps. Lectures on Choquet's theorem. Second edition. Lecture Notes in Mathematics, 1757. Springer-Verlag, Berlin, 2001. viii+124 pp.

\bibitem{Ricker81}
Werner Ricker. A new class of convex bodies.
Papers in algebra, analysis and statistics (Hobart, 1981), pp. 333--340, Contemp. Math., 9, Amer. Math. Soc., Providence, R.I., 1981.

\bibitem{Sallee74}
George T.~Sallee. On the indecomposability of the cone. J.~London Math.~Soc 9 (1974), 363--367.


\bibitem{Schneider1985}
Rolf Schneider. “On the Aleksandrov--Fenchel inequality”. In: Discrete ge-
ometry and convexity (New York, 1982). Vol. 440. Ann. New York Acad.
Sci. New York Acad. Sci., New York, 1985, pp. 132–141.

\bibitem{Schneider1988}
Rolf Schneider. On the Aleksandrov--Fenchel in\-equal\-ity invol\-ving zono\-ids. Geom. Dedicata 27 (1988), no. 1, 113--126.

\bibitem{Schneider1996}
Rolf Schneider. Simple valuations on convex bodies. Mathematika 43 (1996), 32--39.

\bibitem{Schneider}
Rolf Schneider. Convex Bodies: the Brunn--Minkowski Theory. Second expanded edition.
Vol. 151. Encyclopedia of Mathematics and its Applications. Cambridge
University Press, Cambridge, 2014, pp. xxii+736.

\bibitem{SW}
Rolf Schneider and Wolfgang Weil.  Stochastic and integral geometry. Probability and its Applications (New York). Springer-Verlag, Berlin, 2008. xii+693 pp.

\bibitem{SvH19}
Yair Shenfeld and Ramon van Handel. Mixed volumes and the Bochner method. Proc. Amer. Math. Soc. 147 (2019), no. 12, 5385--5402.


\bibitem{SvH22}
Yair Shenfeld and Ramon van Handel. The extremals of Minkowski's quadratic inequality.
Duke Math. J. 171 (2022), no. 4, 957--1027.


\bibitem{SvH23+}
Yair Shenfeld and Ramon van Handel. The Extremals of the Alex\-an\-drov--Fenchel inequality for convex polytopes. ArXiv:2011.04059. (To appear in Acta Mathematica.)



\bibitem{Varadarajan1958}
Veeravalli S.~Varadarajan. On the convergence of sample probability distributions.  Sankhya 19 (1958), 23--26.

\bibitem{Wa18}
Yu Wang. A remark on the Alexandrov--Fenchel inequality.
J. Funct. Anal. 274 (2018), no. 7, 2061--2088.

\bibitem{Weil76}
Wolfgang Weil. Kontinuierliche Linearkombinationen von Strecken. Math. Z. 148 (1976), 71--84.



\end{thebibliography}
\end{document}